\documentclass{article}
\usepackage[utf8]{inputenc}

\usepackage{algorithm2e}
\usepackage{xcolor}

\SetAlFnt{\footnotesize}

\SetAlCapSty{xAlCapSty}

\SetCommentSty{xCommentSty}

\SetNlSty{mynlfont}{}{} 

\LinesNumbered

\SetSideCommentRight

\DontPrintSemicolon

\RestyleAlgo{algoruled}

\usepackage[driverfallback=dvipdfm,
bookmarks=true,    
bookmarksopen=true,
bookmarksopenlevel=\maxdimen,
bookmarksdepth=10
]{hyperref}
\hypersetup{
pdfcenterwindow=true, 
pdfpagelabels=true,
pagebackref=true,            
hypertexnames=true,
plainpages=false,
unicode=false,     
pdftoolbar=true,                
    pdfmenubar=true,          
    pdffitwindow=false,         
    pdfstartview={FitH},        
    pdftitle={Accuracy- Source Term Quadrature for 
 Third-Order Edge-Based Discretization},
    pdfauthor={Yi Liu and Hiroaki Nishikawa},    
    pdfsubject={Accuracy-Preserving Source Term Quadrature for 
 Third-Order Edge-Based Discretization},       
    pdfcreator={Hiroaki Nishikawa},   
    pdfproducer={Hiroaki Nishikawa}, 
    baseurl={http://www.hiroakinishikawa.com},
    pdfkeywords={}, 
    pdfnewwindow=true,     
    colorlinks=true,       
    linkcolor=black,       
    citecolor=black,       
    filecolor=black,       
    urlcolor=black,        
  breaklinks=true,
  hyperfigures=true,
  backref=true,
breaklinks=false, 
pdfpagelabels,      
pagebackref,        
hypertexnames=true, 
plainpages=false,   
naturalnames        
}
\usepackage[all]{hypcap}

\usepackage[thinlines]{easytable}

\usepackage{color}      
\usepackage{placeins}
\usepackage{indentfirst}
\usepackage{multirow}

\usepackage{amssymb}
\usepackage{amsmath,latexsym}
\usepackage{url}
\usepackage{enumerate}
\usepackage{graphicx}

\usepackage{booktabs}

\usepackage{subcaption}
\usepackage{cleveref}
\usepackage{mathtools}

\def\dyad{\! \otimes \!}

\usepackage{blindtext,graphicx}
\usepackage[absolute]{textpos}
\begin{document}

\title{\bf \color{black} On False Accuracy Verification of UMUSCL Scheme}

 \author{
{Hiroaki Nishikawa}\thanks{Associate Research Fellow ({hiro@nianet.org}),
 100 Exploration Way, Hampton, VA 23666 USA.}\\
  {\normalsize\itshape National Institute of Aerospace,
 Hampton, VA 23666, USA} 
}

\date{\today}
\maketitle

\begin{abstract} 
In this paper, we reveal a mechanism behind a false accuracy verification encountered with unstructured-grid schemes based on solution reconstruction 
such as UMUSCL. Third- (or higher-) order of accuracy has been reported for the Euler equations in the literature, but UMUSCL is actually second-order 
accurate at best for nonlinear equations. False high-order convergence occurs generally for a scheme that is high order for linear equations but second-order for nonlinear equations. It is caused by unexpected linearization of a target nonlinear equation due to too small of a perturbation added to an exact solution used for accuracy verification. To clarify the mechanism, we begin with a proof that the UMUSCL scheme is third-order accurate only for linear equations. Then, we derive a condition under which the third-order truncation error dominates the second-order error and demonstrate it numerically for Burgers' equation. Similar results are shown for the Euler equations, which disprove some accuracy verification results in the literature. To be genuinely third-order, UMUSCL must be implemented with flux reconstruction. 
\end{abstract}

\section{Introduction}
\label{introduction}

This paper is a sequel to the two previous papers \cite{Nishikawa_3rdMUSCL:2020,Nishikawa_3rdQUICK:2020}, where we clarified the MUSCL and QUICK schemes towards the clarification of economical high-order unstructured-grid schemes for practical computational fluid dynamics (CFD) solvers, e.g., third-order UMUSCL with $\kappa=1/2$ \cite{burg_umuscl:AIAA2005-4999}, $\kappa=1/3$ \cite{VLeer_Ultimate_III:JCP1977,VAN_LEER_MUSCL_AERODYNAMIC:J1985}, or $\kappa=0$ \cite{katz_work:JCP2015,nishikawa_liu_source_quadrature:jcp2017}. 
In this paper, we will clarify one more confusion: the {\color{black} false}  accuracy verification of the UMUSCL scheme. 



The UMUSCL scheme of Burg \cite{burg_umuscl:AIAA2005-4999} is
 generally considered as an unstructured-grid extension 
of Van Leer's $\kappa$-reconstruction scheme \cite{VLeer_Ultimate_III:JCP1977,VAN_LEER_MUSCL_AERODYNAMIC:J1985} and has been widely employed in practical CFD solvers \cite{Burg_etal:AIAA2003-3983,MurayamaYamamoto2008,yang_harris:AIAAJ2016,fun3d_website,JAFM2017:UMUSCL,yang_harris:CCP2018,DementRuffin:aiaa2018-1305,Barakos:IJNMF2018,Ito_etal:JoA2018,WhiteNishikawaBaurle_aiaa2019-0127,fun3d_manual:NASATM2019,XFengMHD2020,scFLOW:Aviation2020} 
with a confusion over the value of $\kappa$ for giving third-order accuracy on regular or one-dimensional grids: $\kappa=1/2$ \cite{burg_umuscl:AIAA2005-4999}, $\kappa=1/3$ \cite{VLeer_Ultimate_III:JCP1977,VAN_LEER_MUSCL_AERODYNAMIC:J1985}, or  
$\kappa=-1/6$ \cite{yang_harris:AIAAJ2016,yang_harris:CCP2018}. 
The confusion arises mainly from taking different combinations of numerical solution and discretization types {\color{black} : numerical solutions stored as point-valued or cell-averaged
solutions; discretizations of a differential form of a conservation law at a point (finite-difference) or of an integral form over a cell (finite-volume)}. 
For example, Burg originally proposed the UMUSCL scheme as a finite-volume scheme with point-valued numerical solutions stored at nodes on an unstructured grid \cite{burg_umuscl:AIAA2005-4999}. {\color{black} Note that we know he used point-valued numerical solutions because that is the only way he could obtain third-order accuracy
with $\kappa=1/2$ for a one-dimensional nonlinear system (see Refs.\cite{Nishikawa_3rdMUSCL:2020,Nishikawa_3rdQUICK:2020} for details).} This is already confusing because the MUSCL scheme is based on cell-averaged numerical solutions, not point-valued solutions (see Ref.\cite{Nishikawa_3rdMUSCL:2020}). As clarified in the previous paper  \cite{Nishikawa_3rdQUICK:2020}, a third-order finite-volume scheme with point-valued numerical solutions is nothing but the QUICK scheme and therefore the UMUSCL scheme should have been called the UQUICK scheme. In fact, third-order accuracy with $\kappa=1/2$, which is true for the QUICK scheme, has been confirmed for a one-dimensional steady conservation law by Burg \cite{burg_umuscl:AIAA2005-4999} (see also Ref.\cite{Nishikawa_3rdQUICK:2020}). Third-order accuracy demonstrated by Burg is genuine but only for one-dimensional problems; it cannot be third-order in multi-dimensions even for Cartesian grids unless the flux is integrated over a face by a high-order quadrature formula. 

 To be even more confusing, in many or perhaps all practical unstructured-grid codes, the UMUSCL scheme is implemented 
 {\color{black} not as a finite-volume scheme but as a point-wise scheme with} the time derivative and source/forcing terms evaluated at a solution point (a node or a cell center)  \cite{fun3d_website,yang_harris:AIAAJ2016,yang_harris:CCP2018,DementRuffin:aiaa2018-1305,Barakos:IJNMF2018,fun3d_manual:NASATM2019}.
 {\color{black} In this paper, we will focus on this particular implementation. We will call it simply the UMUSCL scheme but it should not be confused with Burg's 
 finite-volume UMUSCL scheme.} 
 Then, the fact that the scheme has been shown to achieve up to fourth-order accuracy with a single flux evaluation per face on Cartesian grids \cite{yang_harris:AIAAJ2016,yang_harris:CCP2018,DementRuffin:aiaa2018-1305,Barakos:IJNMF2018} indicates that the scheme is actually a finite-difference scheme, approximating the differential form of a target equation at a solution point. {\color{black} Note that it does not matter how the discretization is derived; the resulting 
 discretization must be high-order as a finite-difference scheme, not as a finite-volume scheme because the time derivative and source/forcing terms are not integrated with high-order 
 quadrature over a cell.} Therefore, {\color{black} the UMUSCL scheme corresponds to neither the MUSCL scheme nor the QUICK scheme} and does not achieve high-order accuracy in the same way as the MUSCL scheme does. As we will show, the UMUSCL scheme is third-order accurate for linear equations with $\kappa=1/3$, but only second-order accurate when applied to nonlinear conservation laws. 
 This feature is common to conservative finite-difference schemes with a flux evaluated with reconstructed solutions as in MUSCL (e.g., those in Refs.\cite{Dervieux_IJNMF1998,DebiezDerieux:CF2000}). 
 Hence, high-order verification results reported in the literature for the UMUSCL scheme applied to the Euler equations {\color{black} are misleading
and/or misinterpreted}. 
  To be genuinely third- or higher-order accurate, it is necessary to directly reconstruct the flux as pointed out for somewhat similar schemes in Ref.\cite{NLV6_INRIA_report:2008,AbalakinBakhvalovKozubskaya:IJNMF2015}. In other words, the UMUSCL scheme of Refs.\cite{yang_harris:AIAAJ2016,yang_harris:CCP2018,DementRuffin:aiaa2018-1305} can be easily made third- or higher-order by direct flux reconstruction. 
 
{\color{black} It is worth pointing out that a similar unstructured-grid MUSCL scheme had already been proposed in Refs.\cite{Dervieux_IJNMF1998,DebiezDerieux:CF2000,CamarriSalvettiKoobusDerieux:CF2004}. The scheme involves a parameter $\beta$ (instead of $\kappa$) and is called the 
$\beta$-scheme. Its basic form is equivalent to the UMUSCL scheme with $\kappa=1-2\beta$. In fact, the $\beta$-scheme is proposed as a finite-difference-like scheme 
with point-valued solutions at nodes, not as a finite-volume scheme, and has been shown to achieve third-order accuracy on a regular 
grid with $\beta=1/3$ (i.e., $\kappa=1/3$, not $\kappa=1/2$) for linear equations \cite{DebiezDerieux:CF2000}. Therefore, all discussions in this paper will equally apply to 
the $\beta$-scheme.  
}
 
It should be noted also that the UMUSCL scheme has been shown to bring significant improvements to complex flow simulations and thus it is indeed useful   \cite{yang_harris:AIAAJ2016,yang_harris:CCP2018,DementRuffin:aiaa2018-1305}. It is just that the improved resolution is largely due to reduced dissipation 
by high-order solution reconstruction, not by high-order accuracy (see Refs.\cite{Dervieux_IJNMF1998,DebiezDerieux:CF2000,CamarriSalvettiKoobusDerieux:CF2004} for a relevant discussion). For scale-resolving turbulent-flow simulations, however, it would be strongly desired to achieve high-order accuracy since a high-order scheme is expected to be much more efficient than second-order methods on highly refined grid. 
It is therefore important to {\color{black} reveal the accuracy limitation of UMUSCL and develop a genuinely high-order version.}

To accomplish the task, we first show that the UMUSCL scheme, {\color{black} in the form typically implemented in a practical code with the time derivative and source/forcing 
terms evaluated at a solution point,} is equivalent to the high-order conservative scheme of Shu and Osher \cite{Shu_Osher_Efficient_ENO_II_JCP1989} only for linear equations and thus cannot be high-order for nonlinear equations. Then, we reveal the mechanism behind the {\color{black} false}  third-order error convergence and derive a condition under which a third-order truncation error dominates the second-order one. Finally, we numerically demonstrate the {\color{black} false}  third-order accuracy of UMUSCL and genuine third-order accuracy of a flux-reconstruction version for the Burgers and Euler equations. Our focus is on third-order accuracy for smooth solutions, which is sufficient to illustrate the problem. A monotonicity property (which is important for discontinuous solutions) and higher-order accuracy are beyond the scope of the paper and will be discussed elsewhere. {\color{black} Accuracy verification will be performed with exact solutions or by the method of 
manufactured solutions \cite{OberkampfRoy_2010}. See Ref.\cite{OberkampfRoy_2010} for a deeper background on verification; but the issue discussed in this paper is not 
addressed in Ref.\cite{OberkampfRoy_2010}.}
 
The paper is organized as follows. In Section 2, we will explain why the UMUSCL scheme cannot be third-order accurate for nonlinear
equations unless the flux is directly reconstructed. 
In Section 3, we will describe the mechanism of how a target nonlinear equation is linearized by an exact solution. 
In Section 4, we will present numerical results to confirm the {\color{black} false}  third-order error convergence of UMUSCL
and genuine third-order accuracy with flux reconstruction.  
In Section 5, we conclude the paper with remarks.

\section{Second/Third-Order for Linear/Nonlinear Equations}
\label{second_third}

\subsection{UMUSCL is second-order accurate}
\label{umuscl_schemes} 

Consider a conservation law in two dimensions:
\begin{eqnarray}
   \partial_t {u} +  \partial_x {f}  + \partial_y g = {s}(x,y),
   \label{burgers_twod}
\end{eqnarray}
where $u$ is a solution variable, $f$ and $g$ are fluxes, and $s$ is a forcing function.
Our target scheme is the UMUSCL scheme with point-valued numerical solutions stored 
at nodes and point evaluations of the time derivative and forcing terms as widely used in practical unstructured-grid 
solvers \cite{yang_harris:AIAAJ2016,Burg_etal:AIAA2003-3983,fun3d_website,yang_harris:CCP2018,DementRuffin:aiaa2018-1305,Barakos:IJNMF2018}. For a general unstructured
grid, the UMUSCL scheme is given at a node $j$ by 
\begin{eqnarray}
\frac{d u_{j} }{dt}   + \frac{1}{V_j} \sum_{ k \in \{ k_j \}} \phi_{jk}(u_L,u_R)  A_{jk} = s_{j},
\label{semi_discrete_general}
\end{eqnarray} 
where $u_j$ is a point-valued solution at the node $j$, $V_j$ is the measure of the {\color{black} dual} control volume around the node $j$, $A_{jk}$ is the area (length in two dimensions) of the face between $j$ and $k$,  $\{ k_j \}$ is a set of neighbor nodes of $j$, $\phi_{jk}$ is a numerical flux along a face normal 
as a function of two states $u_L$ and $u_R$ reconstructed from the nodes $j$ and $k$, respectively, by the so-called UMUSCL reconstruction scheme \cite{burg_umuscl:AIAA2005-4999}, and $s_j$ is a point evaluation of $s(x,y)$ at $j$. {\color{black} Note that we store numerical solutions as point values at nodes (there are no cell-averaged 
solutions stored anywhere) and the time derivative $ \partial_t {u}$ is simply evaluated at the node as $\frac{d u_{j} }{dt} $. }

As clarified in the previous paper \cite{Nishikawa_3rdQUICK:2020} and will be discussed in detail later, the UMUSCL scheme corresponds to {\color{black} a generalized version of} the QUICKEST scheme of Leonard  \cite{Leonard_QUICK_CMAME1979} as well as the $\kappa$-family finite-difference scheme of Van Leer \cite{VLeer_Ultimate_III:JCP1977} or more generally the conservative finite-difference scheme of Shu and Osher \cite{Shu_Osher_Efficient_ENO_II_JCP1989} in one dimension or on Cartesian grids. 
A major advantage of the UMUSCL scheme implemented in the form (\ref{semi_discrete_general}) is its efficiency: it can achieve third-order accuracy on regular grids (at least for linear equations) whereas the original finite-volume version \cite{burg_umuscl:AIAA2005-4999} requires high-order flux quadrature to achieve third-order accuracy even on regular grids and also requires a consistent treatment of the cell-averaged time derivative for unsteady problems \cite{Nishikawa_3rdQUICK:2020}. A further discussion on efficiency will be given in a subsequent paper.
 
\begin{figure}[t]
    \centering 
        \includegraphics[width=0.48\textwidth]{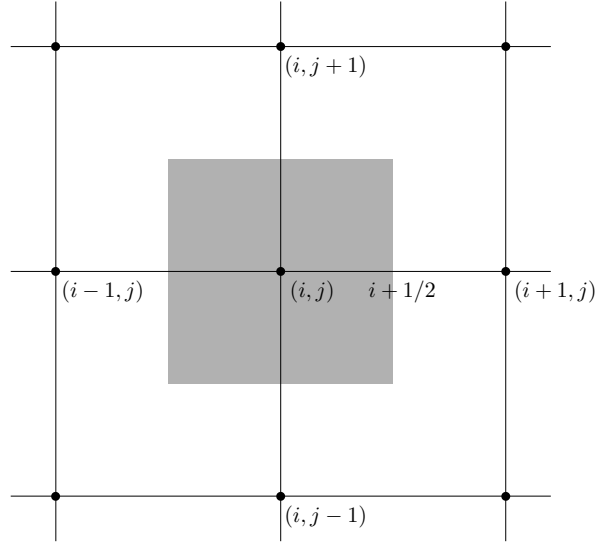} 
            \caption{
\label{fig:oned_fv_data_quad_inter}%
A control volume around a node $(i,j)$, which can also be interpreted as a cell center, on a Cartesian grid. 
} 
\end{figure}

In this paper, we focus on a Cartesian grid as shown in Figure \ref{fig:oned_fv_data_quad_inter}, which is sufficient to illustrate the {\color{black} false}  accuracy problem, and denote the numerical solution at a node $(i,j)$ by $u_{i,j}$. Then, the UMUSCL scheme reduces to
\begin{eqnarray}
\frac{d u_{i,j}}{dt }  + \frac{ F_{i+1/2,j} -   F_{i-1/2,j}  }{h} +  \frac{  G_{i,j+1/2} -  G_{i,j-1/2}    }{h}= s_{i,j},
\label{semi_discrete}
\end{eqnarray} 
where $x_{i+1,j} - x_{i,j} = y_{i,j+1} - y_{i,j} = h$, $s_{i,j} = s(x_{i,j} , y_{i,j} )$, and $F$ and $G$ are numerical fluxes defined at a face,
\begin{eqnarray}
F(u_L,u_R) =  \frac{1}{2} \left[  f(u_L)    + f(u_R)     \right]  - \frac{D}{2}  ( u_R - u_L), 
\label{umuscl_numerical_flux}
\end{eqnarray}
with the dissipation coefficient $D= |\partial f/ \partial u|$, and $u_L$ and $u_R$ are solutions reconstructed with Van Leer's $\kappa$-reconstruction scheme \cite{VLeer_Ultimate_III:JCP1977,VAN_LEER_MUSCL_AERODYNAMIC:J1985} implemented in the form applicable to general unstructured grids 
as suggested by Burg \cite{burg_umuscl:AIAA2005-4999}: e.g., at a face $(i+1/2,j)$, 
\begin{eqnarray}
{u}_L &=& \kappa  \frac{ {u}_{i} +  {u}_{i+1}}{2} +   (1-\kappa) \left[  {u}_j  + \frac{h}{2}(u_x)_{i}  \right] , \quad  (u_x)_{i}  = \frac{   {u}_{i+1} -  {u}_{i-1}  }{2h}, 
\label{umuscl_L2}\\ [1ex]
{u}_R &=& \kappa  \frac{ {u}_{i+1} +  {u}_{i}}{2} +   (1-\kappa) \left[  {u}_{i+1} - \frac{h}{2}(u_x)_{i+1}  \right], \quad (u_x)_{i+1}  = \frac{   {u}_{i+2} -  {u}_{i}  }{2h},
 \label{umuscl_R2}
\end{eqnarray}
where the other subscript $j$ has been omitted for brevity and $\kappa$ is a parameter. 
The other numerical flux $G$ is computed in a similar manner in the $y$-direction. Refs.\cite{yang_harris:AIAAJ2016,yang_harris:CCP2018,DementRuffin:aiaa2018-1305,Barakos:IJNMF2018} add extra terms to the above reconstruction schemes in an attempt to achieve even higher-order accuracy. 
As we will show, it is the flux evaluation with the reconstructed solution $f(u_L)$, not the reconstruction scheme itself, that leads to accuracy deterioration for nonlinear equations. Therefore, such extra terms are irrelevant to the discussion here and thus not considered. 

At first glance, the UMUSCL scheme (\ref{semi_discrete}) looks like a finite-difference scheme approximating the differential form (\ref{burgers_twod}) at the node $(i,j)$. 
As such, it is second-order accurate for any $\kappa$ because it consists of the central difference approximations: 
\begin{eqnarray}
 \frac{ F_{i+1/2,j} -   F_{i-1/2,j}  }{h}  \approx   \partial_x {f}   + O(h^2), \quad
  \frac{  G_{i,j+1/2} -  G_{i,j-1/2}    }{h}  \approx   \partial_y {g}   + O(h^2),
  \label{FD_second_order_at_best}
\end{eqnarray}
no matter how accurate the fluxes $F$ and $G$ are. 
 However, as mentioned in Introduction, third- or higher-order accuracy has been observed with $\kappa=1/3$ as
reported in the literature  \cite{yang_harris:AIAAJ2016,yang_harris:CCP2018,DementRuffin:aiaa2018-1305}. To clarify the apparent contradiction, we 
will show that the scheme (\ref{semi_discrete}) has a slightly different interpretation and does achieve third-order accuracy for linear conservation laws, at least. 




 
\subsection{UMUSCL can be third-order accurate for linear equations}
\label{umuscl_3rd_for_linear}

In 1989, Shu and Osher \cite{Shu_Osher_Efficient_ENO_II_JCP1989} proposed a high-order conservative finite-difference method based on the following 
relationship (which is exact on uniform grids; see Ref.\cite{Merryman:JSC2003} for a discussion on its accuracy on non-uniform grids):
\begin{eqnarray}
    \partial_x {f}   =   \frac{ \tilde{F}_{i+1/2,j} -   \tilde{F}_{i-1/2,j}  }{h},
    \label{ShuOsher1989}
\end{eqnarray}
where $\tilde{F}$ is a function whose cell-average in the $x$-direction is $f(u_{i,j})$, which can be easily derived as 
\begin{eqnarray}
    \partial_x {f}   =  \partial_x \left(  \frac{1}{h} \int_{x-h/2}^{x+h/2}  \tilde{F} dx \right)  =  \frac{ \tilde{F}_{i+1/2,j} -   \tilde{F}_{i-1/2,j}  }{h}.
\end{eqnarray}
Therefore, a high-order scheme can be constructed in the form: 
\begin{eqnarray}
\frac{d u_{i,j}}{dt }  + \frac{ \tilde{F}_{i+1/2,j} -  \tilde{F}_{i-1/2,j}  }{h} +  \frac{  \tilde{G}_{i,j+1/2} -  \tilde{G}_{i,j-1/2}    }{h}= s_{i,j},
\label{semi_discrete_02}
\end{eqnarray} 
by accurately reconstructing $\tilde{F}$ and $\tilde{G}$. Then, since these fluxes are point values at faces and the data available at nodes are their cell-averages, 
we have to reconstruct the point-valued flux, e.g., $\tilde{F}$, from its cell-averages $f_{i,j} = f(u_{i,j})$. 
That is what Van Leer's $\kappa$-reconstruction scheme does exactly for a quadratic function with $\kappa=1/3$ \cite{Nishikawa_3rdMUSCL:2020,VLeer_Ultimate_III:JCP1977,VAN_LEER_MUSCL_AERODYNAMIC:J1985}: e.g., again at a face $(i+1/2,j)$,   
  \begin{eqnarray}
{f}_L &=& \kappa \,  \frac{ {f}_{i} +  {f}_{i+1}}{2} +   (1-\kappa) \left[  {f}_j  + \frac{h}{2}(f_x)_{i}  \right] , \quad  (f_x)_{i}  = \frac{   {f}_{i+1} -  {f}_{i-1}  }{2h}, 
\label{flux_reconstruct_fL} \\ [1ex]
{f}_R &=& \kappa \, \frac{ {f}_{i+1} +  {f}_{i}}{2} +   (1-\kappa) \left[  {f}_{i+1} - \frac{h}{2}(f_x)_{i+1}  \right], \quad (f_x)_{i+1}  = \frac{   {f}_{i+2} -  {f}_{i}  }{2h},
\label{flux_reconstruct_fR}
 \label{umuscl_L2_flux}
\end{eqnarray}
where we have denoted, for brevity, ${f}_{i-1}  = f(u_{i-1,j})$, ${f}_{i}  = f(u_{i,j})$, ${f}_{i+1}  = f(u_{i+1,j})$, ${f}_{i+2}  = f(u_{i+2,j})$.
{\color{black} Note that the above expressions may not look like but are equivalent to the original $\kappa$-scheme (the terms 
proportional to $\kappa$ form a quadratic term; see Ref.\cite{Nishikawa_3rdMUSCL:2020}). }
Then, with the following flux function, 
\begin{eqnarray}
\tilde{F}(f_L, f_R, u_L,u_R) =  \frac{1}{2} \left[  f_L    + f_R   \right]  - \frac{D}{2}  ( u_R - u_L),
\label{cfd_numerical_flux}
\end{eqnarray}
{\color{black} where $D= |\partial f/ \partial u|$,} the conservative finite-difference scheme (\ref{semi_discrete_02}) achieves third-order accuracy with $\kappa=1/3$, which has been verified analytically and numerically in the previous paper \cite{Nishikawa_3rdQUICK:2020}.


Notice now that the scheme (\ref{semi_discrete_02}) is equivalent to the UMUSCL scheme, 
{\color{black} which evaluates the left and right fluxes with reconstructed solutions,  $f(u_L)$ and $ f(u_R)$,}
 for a linear equation, i.e.,  $f = a u$ with a constant $a$, 
because we then have 
\begin{eqnarray}
  f(u_L) = f_L , \quad   f(u_R) = f_R ,
\end{eqnarray}
and thus
\begin{eqnarray}
 \tilde{F}(f_L, f_R, u_L,u_R) =  F(u_L,u_R).
\end{eqnarray}
That is, the solution reconstruction is equivalent to the flux reconstruction for a linear flux. 
Therefore, the UMUSCL scheme achieves third-order accuracy with $\kappa=1/3$ for a linear equation on Cartesian grids.

\subsection{UMUSCL cannot be third-order accurate for nonlinear equations}
\label{umuscl_2nd_for_nonlinear} 

For nonlinear fluxes,  the solution reconstruction is no longer equivalent to the flux reconstruction:
\begin{eqnarray}
f(u_L) \ne   f_L  , \quad  f(u_R)  \ne f_R,
\end{eqnarray}
and thus
\begin{eqnarray}
 \tilde{F}(f_L, f_R, u_L,u_R) \ne  F(u_L,u_R).
\end{eqnarray}
Therefore, the UMUSCL scheme is not equivalent to the conservative finite-difference scheme and simply reduces to 
a second-order finite-difference scheme as in Equation (\ref{FD_second_order_at_best}). 
It implies also that any point-wise numerical scheme with $F(u_L,u_R)$, e.g., those in Refs.\cite{Dervieux_IJNMF1998,DebiezDerieux:CF2000}, is second-order accurate at best. 

To see this more clearly, consider a quadratic function whose cell-average is the flux at a node $j$: 
\begin{eqnarray}
  \tilde{f}(x)  =   {f}_{i} + (f_x)_{i}  (x - x_{i} ) + \frac{1}{2} (f_{xx})_{i} \left[  (x - x_{i} )^2 -   \frac{h^2}{12} \right]
  , \quad (f_{xx})_{i}  = \frac{   {f}_{i+1} - 2 {f}_{i}  + f_{i-1} }{h^2},
\end{eqnarray}
which gives at a face $i+1/2$ 
\begin{eqnarray}
  \tilde{f}_{i+1/2}  =   {f}_{i} + \frac{h}{2} (f_x)_{i}  + \frac{h^2}{12} (f_{xx})_{i},
  \label{exact_quadratic_flux_at_right}
\end{eqnarray}
which can be expanded, for Burgers' equation $f=u^2/2$, as 
\begin{eqnarray}
  \tilde{f}_{i+1/2}  =  \frac{u_i^2}{2}  + \frac{1}{2} ( u_i \partial_{x} u ) h + \frac{1}{12}   \left[     (\partial_{x} u)^2  + u_i \partial_{xx} u   \right] h^2
  +O(h^3).
  \label{exact_fL}
\end{eqnarray}
This at least needs to be matched exactly by a numerical flux in order to generate a third-order scheme. For the direct flux reconstruction with $\kappa=1/3$, we find
\begin{eqnarray}
f_L &=& 
\kappa  \frac{ {f}_{i} +  {f}_{i+1}}{2} +   (1-\kappa) \left[  {f}_j  + \frac{1}{2}(f_x)_{i}  h \right]   \nonumber   \\ [2ex]
&=& \frac{u_i^2}{2} +   \frac{1}{2}  ( u_i \partial_{x} u ) h +  \frac{1}{12}   \left[     (\partial_{x} u)^2  + u_i \partial_{xx} u   \right] h^2+ O(h^3),
\end{eqnarray}  
which matches the exact expansion (\ref{exact_fL}). 
 However, for the UMUSCL scheme, we find with $\kappa=1/3$, 
\begin{eqnarray}
  f(u_L)  
  &=& \frac{1}{2} 
   \left\{ 
    \kappa  \frac{ {u}_{i} +  {u}_{i+1}}{2} +   (1-\kappa) \left[  {u}_j  + \frac{h}{2}(u_x)_{i}  \right] 
    \right\}^2  \nonumber   \\ [2ex]
&=& \frac{u_i^2}{2} +   \frac{1}{2}  ( u_i \partial_{x} u )h  + 
 \frac{1}{12}   \left[  { \color{black}  \frac{3}{2} } (\partial_{x} u)^2  + u_i \partial_{xx} u   \right]   h^2 + O(h^3),
\end{eqnarray}  
which does not match the exact expansion (\ref{exact_fL}) in the quadratic term. Therefore, the UMUSCL scheme is only second-order accurate for nonlinear equations. To understand it more deeply, one has to understand that the flux at a face needs to be computed as a function whose cell-average is $f_i$ and directly reconstructed from the cell-averages $f_i$, $f_{i+1}$, etc., whereas in the latter the flux is evaluated with the face value of a function whose cell-average is $u_i$, reconstructed from the cell averages $u_i$, $u_{i+1}$, etc.; it is well known that they differ by a second-order error \cite{Margolin:ShockWaves2019} (the error in the above can be predicted by the theorem in Ref.\cite{Margolin:ShockWaves2019}). See also Ref.\cite{Nishikawa_3rdQUICK:2020} for discussions on the QUICKEST scheme, which is equivalent to the UMUSCL scheme. A further discussion will be given in a subsequent paper. 

Nonetheless, it should also be clear by now that the UMUSCL scheme can be made third-order for nonlinear equations by performing the flux reconstruction to replace $f(u_L)$ by $f_L$. 
 In the rest of the paper, the flux-reconstruction version of the UMUSCL scheme with the 
numerical flux (\ref{cfd_numerical_flux}) will be referred to as FSR (flux and solution reconstruction). The word `UMUSCL' is not used here because it is not MUSCL as mentioned earlier. We will come back to this scheme later in numerical experiments. Note that the necessity of flux reconstruction was discussed in Ref.\cite{NLV6_INRIA_report:2008} for a similar reconstruction-based scheme; but they did not rely on the exact relation (\ref{ShuOsher1989}) and instead directly constructed a high-order approximation to the flux derivative. 

{\color{black} It is confusing that} third-order accuracy of the UMUSCL scheme has been observed for the nonlinear Euler equations 
 and such accuracy verification results have been used as a confirmation of third- and higher-order accuracy of the 
UMUSCL scheme \cite{yang_harris:AIAAJ2016,yang_harris:CCP2018,Nishikawa_FANG_AQ:Aviation2020}. 
Note that in Ref.\cite{Nishikawa_FANG_AQ:Aviation2020}, the author observed third-order accuracy with a similar scheme but later realized that it was {\color{black} false}  for the reason discussed in this paper. Note also that Refs.\cite{yang_harris:AIAAJ2016,yang_harris:CCP2018} mention high-order flux quadrature but it 
{\color{black} is not sufficient}
 because they treat source terms as point values (thus not finite-volume). These results are not genuine and thus misleading (no other accuracy verifications are shown in Refs.\cite{yang_harris:AIAAJ2016,yang_harris:CCP2018,Nishikawa_FANG_AQ:Aviation2020}). 
Here, the problem exists in the accuracy verification process as we will discuss in the next section.

\section{Unexpected Linearization by Exact Solution}
\label{false_accuracy}


\subsection{Burgers equation}
\label{burgers} 

The problem is illustrated for the steady Burgers equation, a representative of nonlinear conservation laws, in one dimension:
\begin{eqnarray}
    \partial_x {f} = {s}(x),
   \label{burgers}
\end{eqnarray}
where $f = u^2/2$ and $s(x)$ is defined such that a chosen function $u_e(x)$ is made an exact solution: 
$s(x) = \partial_x  ( u_e^2 / 2)  $. Consider a typical exact solution used in accuracy verification tests:
\begin{eqnarray}
    u_{e}(x)   = u_\infty +  \epsilon u_p(x), 
    \label{special_exact_solution}
\end{eqnarray}
where $u_\infty$ and $ \epsilon$ are positive constants, and $u_p(x)$ is an arbitrary smooth function (e.g., $u_p(x)=\sin(x)$). 
This is the type of solution that has been used to demonstrate high-order accuracy of the UMUSCL scheme \cite{yang_harris:AIAAJ2016,yang_harris:CCP2018,DementRuffin:aiaa2018-1305,Nishikawa_FANG_AQ:Aviation2020}. 
Consider the UMUSCL scheme for the steady Burgers equation:
\begin{eqnarray}
 \frac{  F_{i+1/2} -   F_{i-1/2}  }{h} = s_{i},
\label{semi_discrete_1d}
\end{eqnarray}
where the fluxes are computed as described in Section \ref{umuscl_schemes}. 
From the discussion in the previous section, we know this scheme cannot be third-order accurate because the Burgers equation is nonlinear. 
However, it is possible that third-order accuracy is observed accidentally when it is tested with the exact solution of the form (\ref{special_exact_solution}) as
we will explain in the next section.

\subsection{Linearized by exact solution}
\label{linearized} 


False third-order error convergence is observed typically when the parameter $ \epsilon$ is small, which effectively 
linearizes the target nonlinear equation and allows a third-order error to dominate. 
This is the mechanism behind the {\color{black} false}  third-order accuracy verification. 
To see this, substitute the exact solution (\ref{special_exact_solution}) into the Burgers equation to get
\begin{eqnarray}
  \epsilon  \partial_x {f'} = {s}(x),
\end{eqnarray}
where 
\begin{eqnarray}
  {f'} =u_\infty u_p +   \epsilon  u_p^2,
\end{eqnarray}
which shows that the leading term of the flux is linear. Therefore, 
the Burgers equation will behave like a linear equation as $\epsilon \rightarrow 0$; and the UMUSCL scheme exhibits third-order accuracy. 
A question arises, then, about how small $\epsilon$ should be in order to observe third-order error convergence, which is the subject for the next section.

\subsection{Dominated by third-order error}
\label{linearized_dominates} 

To estimate how small $\epsilon$ should be for the {\color{black} false}  third-order error convergence to occur, we substitute a smooth exact solution 
into the residual,
\begin{eqnarray}
\quad Res_i =  \frac{  F_{i+1/2} -   F_{i-1/2}  }{h} -  s_{i},
\label{semi_discrete_1d_res}
\end{eqnarray}
and expand it to obtain the truncation error ${\cal TE}$ (see Ref.\cite{Nishikawa_3rdQUICK:2020}): 
\begin{eqnarray}
{\cal TE} = 
 \frac{1}{24} \left[     f_{uuu} (u_x)^3 +  6 \kappa  f_{uu} u_x  u_{xx}   +  2 \left(  3 \kappa - 1   \right) f_u  u_{xxx} \right] h^2    
+ 
\frac{1}{12} \left[  u_i u_{xxxx}  +   u_x u_{xxx}   \right] h^3  + 
 O(h^4). 
\end{eqnarray}
where the second- and third-order terms come from the averaged flux term and the dissipation term in the numerical flux (\ref{umuscl_numerical_flux}), respectively. 
For simplicity, we have assumed $u_i > 0$ for all $i$. 
Observe that the second-order error vanishes with $\kappa=1/3$ for linear equations, for which $f_{uuu}= f_{uu} = 0$.
However, the Burgers equation gives $f_{uu}=1$ and thus the second-order term remains:   
\begin{eqnarray}
{\cal TE} =   \frac{1}{12} ( u_x )( u_{xx}  ) h^2 + 
\frac{1}{12} \left[  u_i u_{xxxx}  +   u_x u_{xxx}   \right] h^3  + O(h^4).
   \label{scheme_burgers_expanded_02}
\end{eqnarray} 
If the smooth exact solution is given by Equation (\ref{special_exact_solution}) with $u_p(x) =  \sin( \omega x)$, where $\omega$ is a constant, then
\begin{eqnarray}
{\cal TE} =   \epsilon \omega^3  \left[ {\cal T}_2     +    {\cal T}_3   \right] h^2 + O(h^4),
\end{eqnarray} 
where
\begin{eqnarray}
  {\cal T}_2 = - \frac{\epsilon \sin( 2 \omega x ) }{24}
, \quad
  {\cal T}_3 =  -\frac{   \omega \left\{  \epsilon ( 1 - 2  \sin^2( \omega x)  )  - u_\infty \sin(\omega x)    \right\}  h  }{12}.
\end{eqnarray} 
Our interest is to see how small $\epsilon$ should be in order for ${\cal T}_3$ to dominate $ {\cal T}_2 $, and thus the scheme 
will become third-order accurate. It suffices then to compare the upper bounds ($\left| {\cal T}_3  \right|$ reaches the maximum at $x=\frac{\pi}{2 \omega}$): 
\begin{eqnarray}
  \left| {\cal T}_2  \right|_{max} =   \frac{\epsilon   }{24}   , 
  \quad
  \left| {\cal T}_3  \right|_{max} =   \frac{ \omega  ( u_\infty + \epsilon  ) h  }{12},
\end{eqnarray} 
from which we find 
\begin{eqnarray}
  \left| {\cal T}_3  \right|_{max} >   \left| {\cal T}_2  \right|_{max}   \quad \mbox{if}  \quad
 \frac{ \epsilon }{u_\infty}< \frac{2 h \omega }{1- 2 h \omega} .
 \label{largest_epsilon_for_3rd}
\end{eqnarray} 
Therefore, the scheme will behave as if it is a third-order scheme when $\epsilon$ satisfies the above condition.
For example, if an accuracy verification study is performed with $\omega=2 \pi$ over a series of grids with the finest grid of 128 nodes in a unit 
domain giving $h =1/127$, we find
\begin{eqnarray}
 \frac{ \epsilon }{u_\infty} < 0.1098136158 \ldots, 
\end{eqnarray} 
which one would easily choose, since the parameter $\epsilon$ is meant to be a perturbation, if not aware of the problem. 
It is important to note that the analysis shows that the order of accuracy changes based on the value of $\epsilon$. This is quite different from a typical observation that 
a lower-order scheme is more accurate on coarse grids, where the order of accuracy does not change. 
{\color{black} Note also that the problem does not occur even on coarse grids if $\epsilon$ is sufficiently large. 
One can estimate a critical value of $h$ by solving $  \left| {\cal T}_3  \right|_{max} >   \left| {\cal T}_2  \right|_{max} $ for $h$:
\begin{eqnarray}
h > \frac{ \epsilon/u_\infty}{ 2 \omega ( \epsilon /u_\infty + 1 )}.
\label{critical_h_estimate}
\end{eqnarray}
For example, for $\omega= 2 \pi$, we find that third-order accuracy will be observed for grids coarser than the one with 
$h = 0.007234315595 \approx 1/127$ for $\epsilon/u_\infty = 0.1$ and 
$ h = 0.02652582384 \approx 1/38$ for $\epsilon/u_\infty = 0.5$. 
In other words, one must use grids finer than $h=1/127$ and $h=1/38$, in the former and latter cases, respectively, in order to avoid false third-order error convergence.
}

\subsection{Euler equations}
\label{linearizez_euler} 

The same mechanism applies to the Euler equations:
\begin{eqnarray}
  \partial_x {\bf f}  = {\bf s}(x) , \quad  {\bf f} = ( \rho, \rho u^2 + p, \rho u H ),
   \label{1D_euler_system}
\end{eqnarray}
where  $\rho$ is the density, $u$ is the velocity, $p$ is the pressure, and $\rho H = \gamma p/ (\gamma-1) + \rho u^2/2 $ with $\gamma=1.4$.
A similar unexpected linearization occurs when an exact solution for the velocity is defined as in Equation (\ref{special_exact_solution}) 
while exact solutions for the density and the pressure can be defined by arbitrary smooth functions.
Substitute the velocity (\ref{special_exact_solution}) into the Euler equations, 
and obtain
\begin{eqnarray}
   \partial_t {\bf u} + \partial_x {\bf f}_1 + \epsilon \partial_x {\bf f}_2  = {\bf s}(x) ,
   \label{1D_NS_Conservative_linearized}
\end{eqnarray}
where
\begin{eqnarray}
 {\bf f}_{1}
 =
 \left[  \begin{array}{c} 
   \rho  u_\infty \\  [1ex]
               \rho u_\infty^2   + p \\ [1ex]
              u_\infty {\cal H} 
              \end{array} \right]
, \quad
   {\bf f}_2 =  
 \left[  \begin{array}{c} 
        \rho   u_p   \\  [1ex]
             2   \rho u_p    u_\infty  +  \epsilon \rho  u_p^2   \\ [1ex]
        u_p   \left[ {\cal H} +     \rho  u_\infty^2   \right] +   \frac{\epsilon \rho}{2} (    3 u_\infty + \epsilon  u_p    )  u_p^2 
              \end{array} \right],
              \quad
              {\cal H} = \frac{\gamma p }{\gamma-1}  + \frac{1}{2} \rho u_\infty^2  .
   \label{split_flux_euler}
\end{eqnarray}
The factor $\epsilon$ in the second flux derivative term indicates that the system will be linear as $\epsilon \rightarrow 0$.
Therefore, we expect to observe third-order error convergence for the Euler equations if the solution (\ref{special_exact_solution})
is used as an exact solution with a small $\epsilon$. For an estimate of the largest $\epsilon$, we do not attempt to derive it  
for the Euler equations and instead demonstrate that the estimate (\ref{largest_epsilon_for_3rd}) serves rather well also for the Euler equations.

Before proceeding, we remark that there are two other cases where the Euler equations are effectively linearized. 
These case are somewhat trivial but have been used for accuracy verification in the literature, e.g., in Ref.\cite{ZhongSheng:CF2020}.
One is a case where the velocity and the pressure are constant (i.e., the so-called entropy wave). In this case, only the density varies in space and the Euler equations 
become linear since the flux $ {\bf f} $ is a linear function of $\rho$. The other is a case where the velocity is constant. 
In this case, only the density and the pressure vary in space and the Euler equations again become linear since the flux $ {\bf f} $ is a linear function of $\rho$ and $p$. 
Therefore, third-order accuracy will be observed with the UMUSCL scheme. In Ref.\cite{ZhongSheng:CF2020}, high-order accuracy of a WENO-type scheme 
based on solution reconstruction is verified for an entropy-wave solution in one dimension. Their scheme seems different from but should be similar to the UMUSCL scheme because the third-order part is equivalent to the UMUSCL scheme with $\kappa=1/3$ and also it is shown to produce high-order solutions on regular grids in two dimensions without high-order flux 
quadrature over each face. Then, it must be a finite-difference scheme and thus flux reconstruction is required for achieve third- and higher-order accuracy.

\section{Numerical Results}
\label{results}
 
In this section, we present numerical results to demonstrate that third-order accuracy is indeed observed for the Burgers equation when the condition (\ref{largest_epsilon_for_3rd}) is satisfied,
that the same is observed for the Euler equations in one dimension, 
and finally that the same is observed in two-dimensional problems. The last part is intended for disproving high-order accuracy reported in the literature.  
In all numerical experiments, the numerical solutions are point values stored at nodes (which can be interpreted as cell centers also). 
  
\subsection{Burgers equation in one dimension}
\label{results_burgers_1d}

We consider the Burgers equation (\ref{burgers}) in a unit domain $x \in [0,1]$ and the exact solution (\ref{special_exact_solution}) with 
\begin{eqnarray}
u_\infty = 0.3 , 
\quad 
u_p(x) = \sin (2  \pi x),
\end{eqnarray}
which corresponds to taking $\omega=2 \pi$, and solve it over a series of grids with $n =$ 16, 32, 64, 128, and {\color{black} 256} nodes. 
To exclude effects of boundaries (which are beyond the scope of this paper), 
we fix the exact solution at two boundary nodes and their neighbors also, and solve the system of nonlinear residual equations
for those stored at the rest of the cells. The dissipation coefficient $D$ is evaluated with the averaged solution at a face: $ (u_L + u_R )/2 $.
An implicit solver based on the exact Jacobian of the first-order scheme is used to solve the residual equations,
which is taken to be converged when the $L_1$ norm of the residual is reduced by nine orders of magnitude from an initial residual norm computed
with the initial solution $u_i = 1$, $n=3,4, \cdots, n-3, n-2$.
See Ref.\cite{nishikawa_liu_jcp2018}, for example, for details of the implicit solver for a one-dimensional finite-volume-type scheme. 

{\color{black} 
To investigate the effect of the parameter $\epsilon$, we consider varying $\epsilon$ as follows,
\begin{eqnarray}
\epsilon  = c_\epsilon u_\infty, \quad c_\epsilon = 0.05,   \, 0.1,  \,  0.3,   \, 0.5,
\label{epsilon_test_formula}
\end{eqnarray}
which correspond, respectively, to the critical mesh spacing (\ref{critical_h_estimate}), denoted by $h_{cr}$, 
 \begin{eqnarray}
h_{cr} =   \, 0.000410  ,    \,  0.00201,   \, 0.00747,  \,    0.0107.
\end{eqnarray} 
Comparing these with the actual grid spacings, 
 \begin{eqnarray}
h  =  0.0666, \, 0.0322 , \, 0.0159, \, 0.00787 , \,    0.00392, 
\end{eqnarray}
we expect that third-order accuracy will be observed for $c_\epsilon = 0.05$ and $0.1$, and second-order accuracy will begin to 
dominate somewhere between the last two grids for $c_\epsilon = 0.3$ and between the third and fourth grids for $c_\epsilon = 0.5$.  
} Also, we consider three values of $\kappa$, $\kappa=0, 1/2$, and $1/3$, to remind that third-order accuracy 
is possible only with $\kappa=1/3$.

\begin{figure}[th!]
    \centering
      \hfill  
                \begin{subfigure}[t]{0.48\textwidth}
        \includegraphics[width=0.9\textwidth]{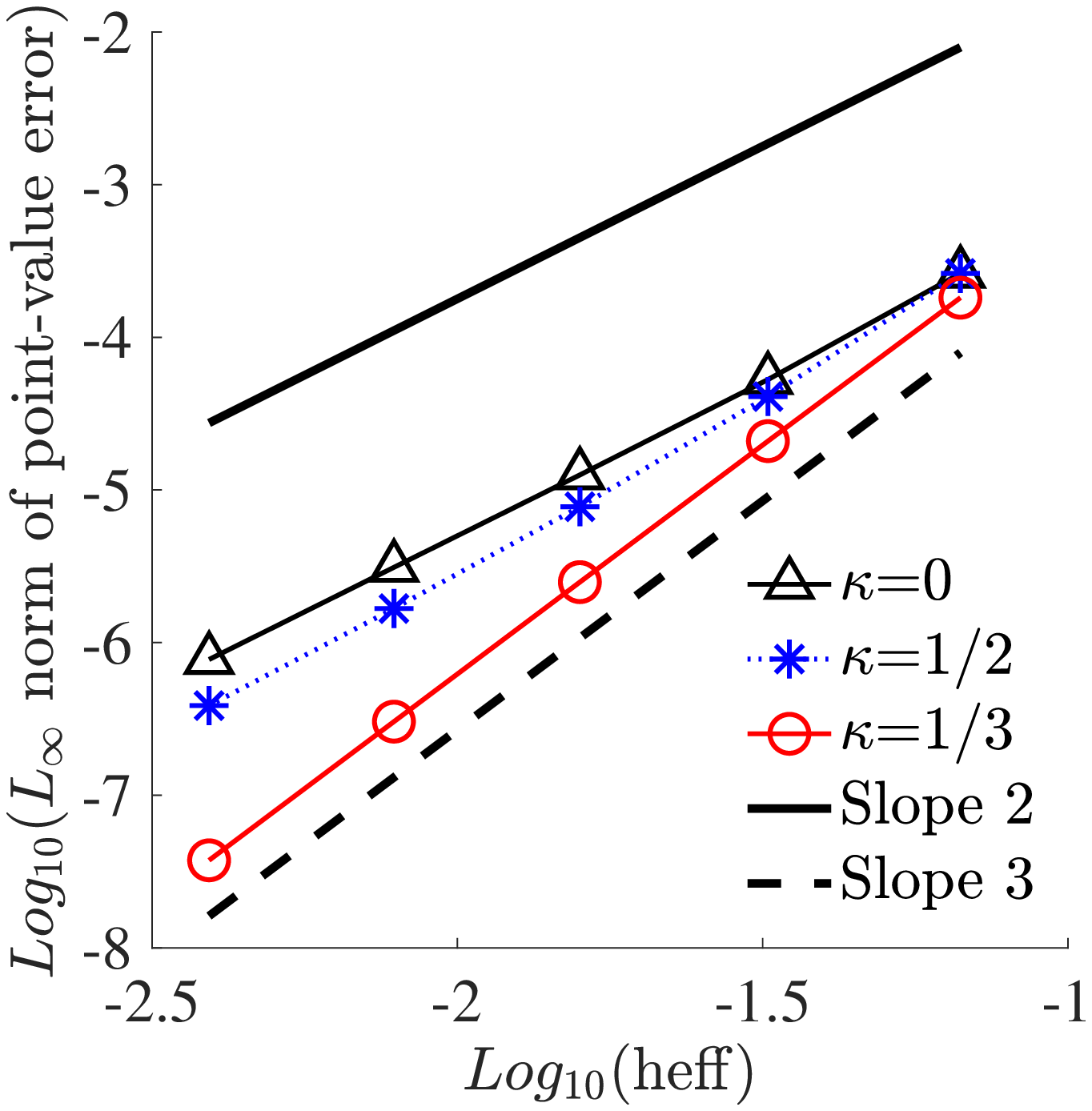}
          \caption{ \color{black} $c_\epsilon = 0.05$.}
       \label{fig:steady_adv_te_c}
      \end{subfigure}
      \hfill
          \begin{subfigure}[t]{0.48\textwidth}
        \includegraphics[width=0.9\textwidth]{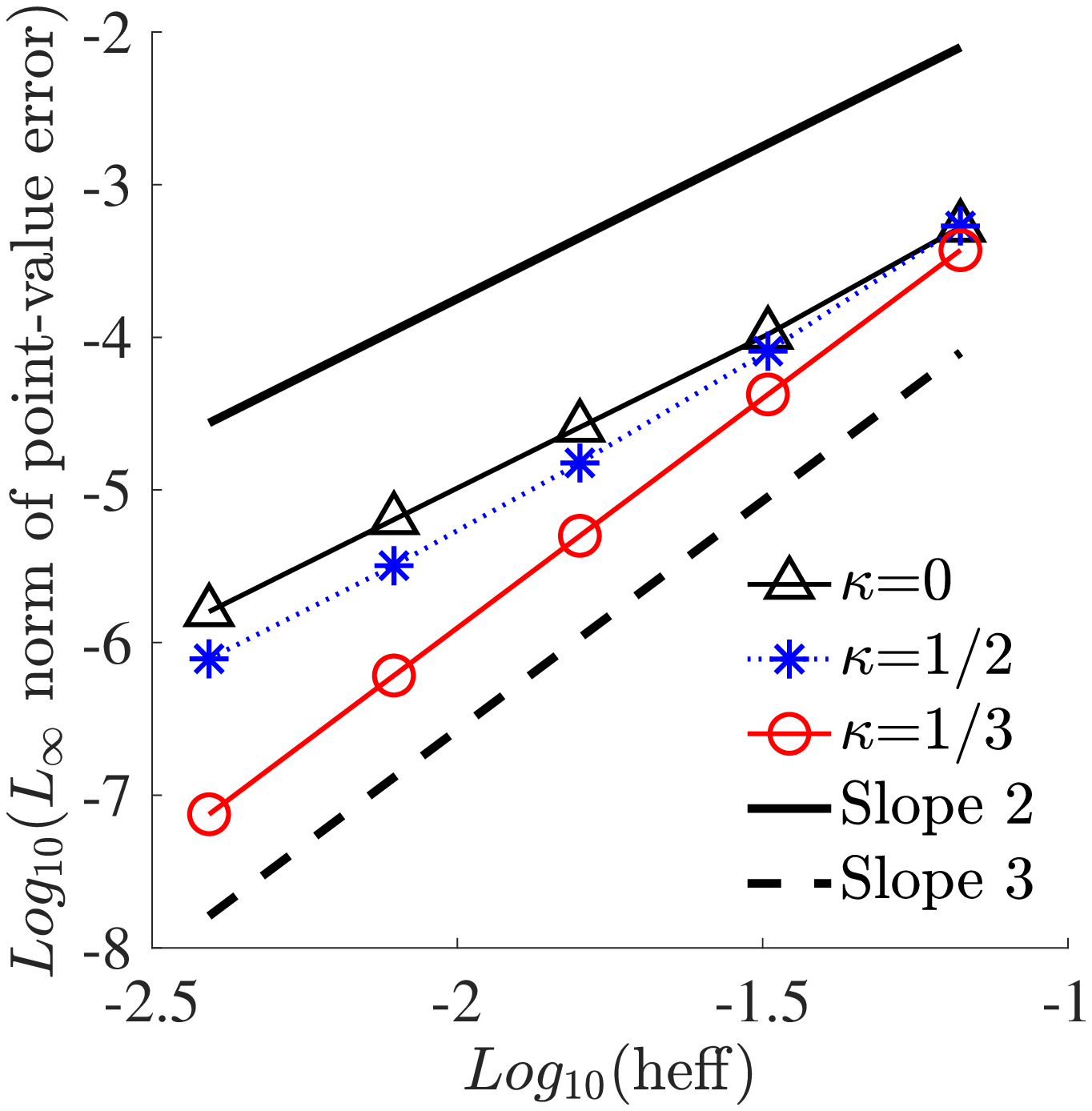}
          \caption{ \color{black} $c_\epsilon = 0.1$.}
       \label{fig:steady_adv_err_ca}
      \end{subfigure}
      \hfill
          \begin{subfigure}[t]{0.48\textwidth}
        \includegraphics[width=0.9\textwidth]{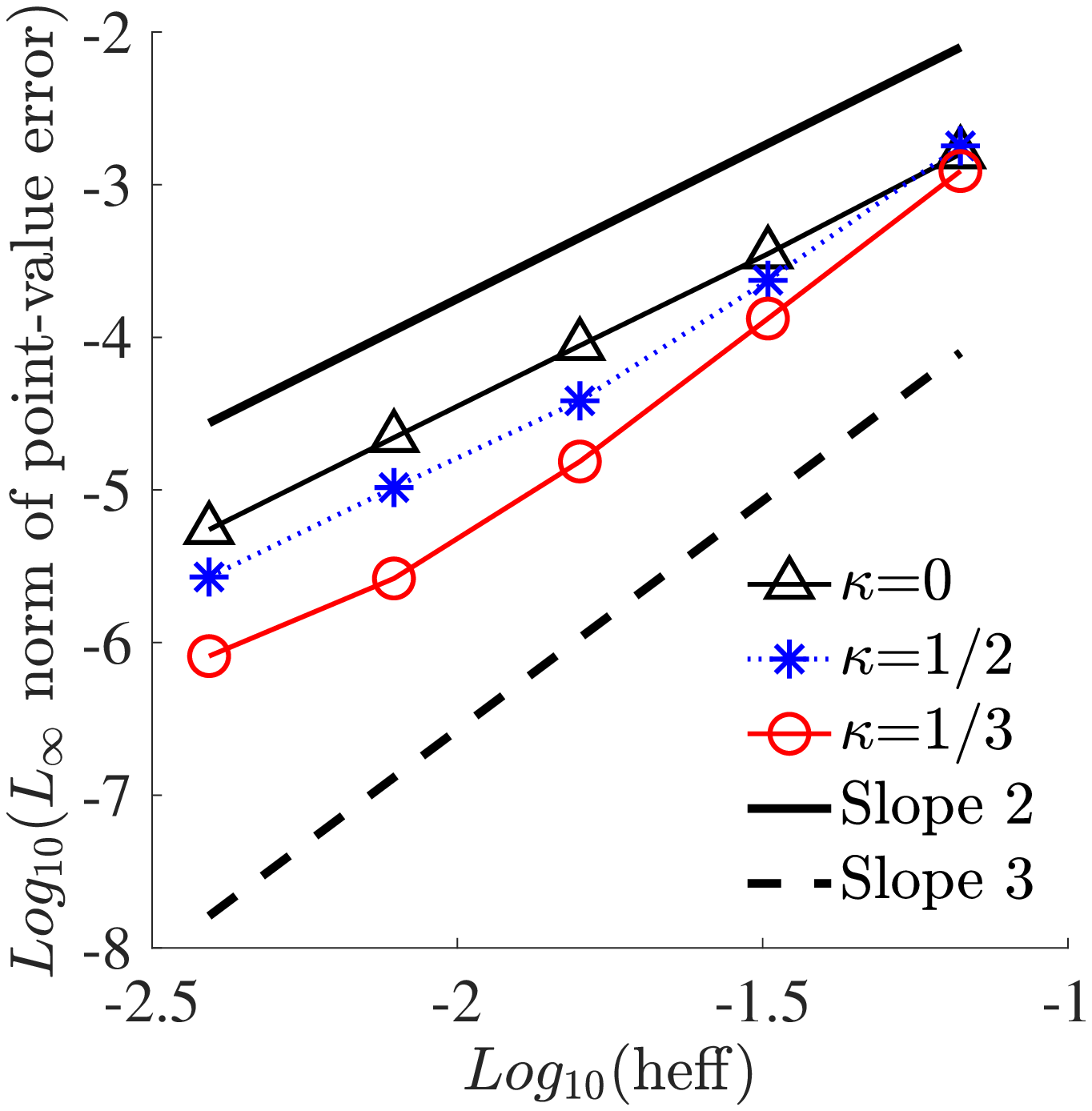}
          \caption{ \color{black} $c_\epsilon = 0.3$.}
       \label{fig:steady_adv_err_ca}
      \end{subfigure}
      \hfill  
                \begin{subfigure}[t]{0.48\textwidth}
        \includegraphics[width=0.9\textwidth]{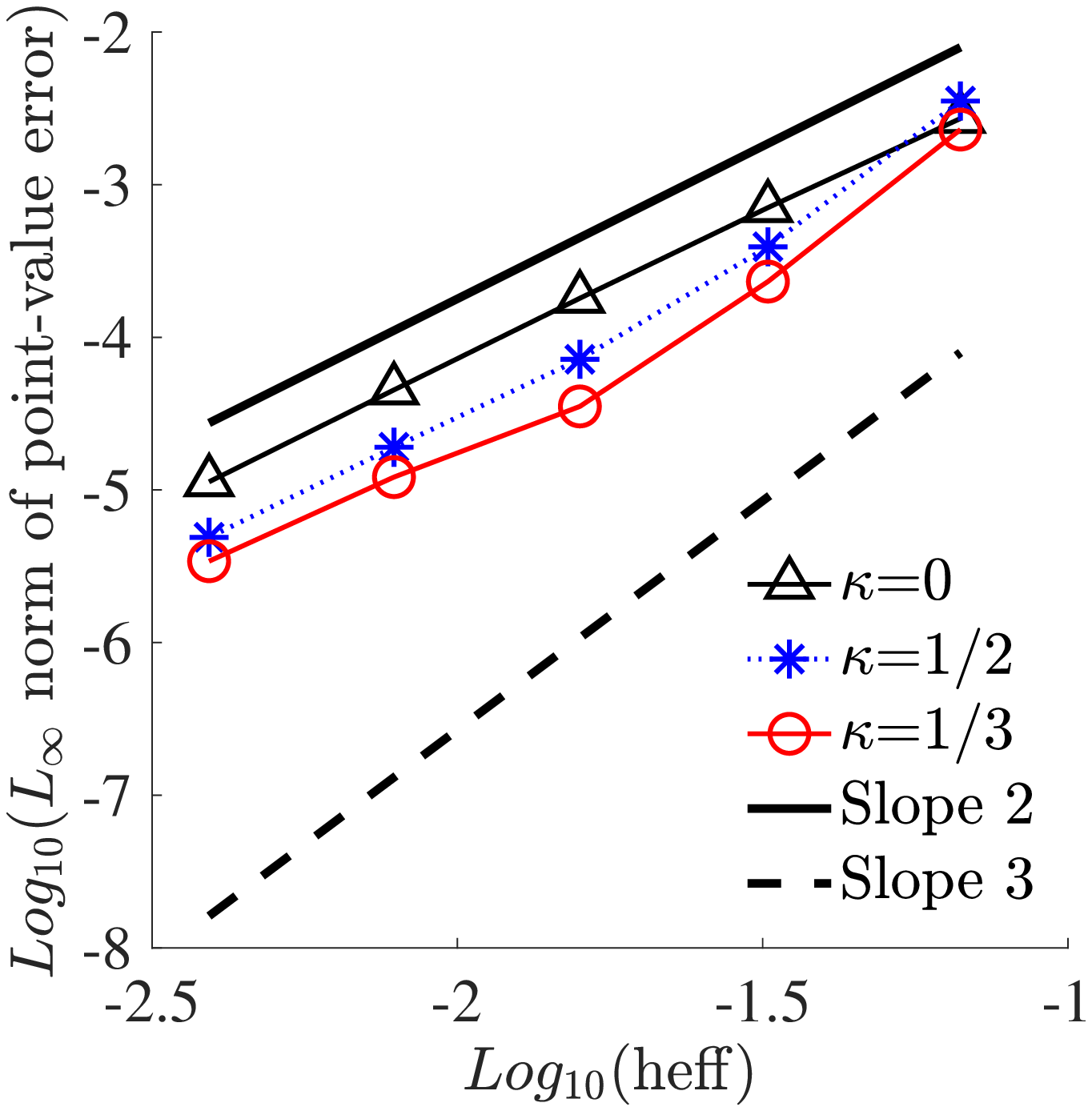}
          \caption{  \color{black}$c_\epsilon = 0.5$.}
       \label{fig:steady_adv_te_p}
      \end{subfigure}
      \hfill    
      \\
            \caption{
\label{fig:steady_adv_error}%
Error convergence for the Burgers equation: 
.
} 
\end{figure}

Error convergence results are shown in Figure \ref{fig:steady_adv_error}, where the discretization error (i.e., solution error) is defined by
\begin{eqnarray}
L_\infty( {\cal E}_p) = \max_{i  \in \{ I \} } | {u}_i - u_i^{exact} |,
\quad
u_i^{exact} = u_e(x_i) = u_\infty + \epsilon \sin (2 \pi x_i ).
\end{eqnarray}
where $\{ I \}$ denotes the set of nodes in a given grid. 
{\color{black} 
As expected, third-order error convergence is observed for $\kappa=1/3$ with $c_\epsilon = 0.05$ and $0.1$, as shown in the first two plots.
Figure \ref{fig:steady_adv_err_ca} shows the result for $c_\epsilon = 0.3$. Second-order convergence is observed in the last two grids as expected.
Finally, Figure \ref{fig:steady_adv_te_p} shows that $c_\epsilon = 0.5$ leads to third-order error convergence begins to deteriorate somewhere
between the third and fourth grids, again as expected.
}

\subsection{Euler equations in one dimension}
\label{results_euler_1d} 

Next, we consider the Euler equations (\ref{1D_euler_system}) with the exact solution defined in terms of the primitive variables:
\begin{eqnarray}
{\bf w}^{exact} (x)
=
\left[
\begin{array}{c}
\rho^{exact} (x) \\
u^{exact} (x)  \\
p^{exact} (x)
\end{array}
\right]
=
\left[
\begin{array}{c}
 1 +  \epsilon_\rho \sin( 2.3 \pi ) \\
u_\infty +  \epsilon  \sin( 2 \pi ) \\
 1 + \epsilon_p \sin( 2.5 \pi x ) 
\end{array}
\right], 
\label{euler_1d_exact_sol}
\end{eqnarray}
and the forcing term defined by ${\bf s}(x) = \partial_x {\bf f}({\bf w}^{exact}(x))$. The amplitudes for the density and the pressure are set $\epsilon_\rho=\epsilon_p=0.2$ unless otherwise stated. For the velocity solution, we set $u_\infty = 0.3 $ and define $\epsilon$ by Equation (\ref{epsilon_test_formula}) to see how well the formula derived for the Burgers equation works for the Euler equations.

In the UMUSCL scheme, we perform the solution reconstruction in the primitive variables ${\bf w}=(\rho,u,p)$ instead of the conservative
variables ${\bf u}=(\rho, \rho u, \rho H - p)$. Then, the residual is given by
\begin{eqnarray}
{\bf Res}_i =  \frac{  {\bf F}_{i+1/2} - {\bf F}_{i-1/2}  }{h} -  {\bf s}_{i},
\label{semi_discrete_euler_1d}
\end{eqnarray}
where the numerical flux is given by the Roe flux \cite{Roe_JCP_1981},
\begin{eqnarray}
{\bf F}({\bf w}_L,{\bf w}_R) =  \frac{1}{2} \left[   {\bf f}({\bf w}_L)   +  {\bf f}({\bf w}_R)   \right]  - \frac{\hat{\bf D}}{2}  \left[  {\bf u}({\bf w}_R) -  {\bf u}({\bf w}_L) \right],
\label{cfd_numerical_flux_euler}
\end{eqnarray}
$\hat{\bf D}$ is the Roe dissipation matrix evaluated with the Roe averages \cite{Roe_JCP_1981}, and the two states ${\bf w}_L$ and ${\bf w}_R$ are computed by
the reconstruction formulas (\ref{umuscl_L2}) and (\ref{umuscl_R2}).  As before, we consider three values of $\kappa$: $\kappa=0, 1/2$, and $1/3$.

The resulting system of steady residual equations is solved by a pseudo-time integration with the three-stage SSP Runge-Kutta scheme \cite{SSP:SIAMReview2001} 
with a local time step at CFL$=0.99$. The solver is taken to be converged when the residual is reduced by seven orders of magnitude in the $L_1$ norm, starting from
the initial norm computed with the initial constant solution: ${\bf w}_i = (1, u_\infty , 1 )$, $i=3,4, \cdots, n-3, n-2$. Again, the exact solutions are specified at $i=1,2 , n-1$, and $n$, to exclude boundary effects.
For simplicity, the discretization error is defined as the maximum of the $L_\infty$ norms among the three variables:
\begin{eqnarray}
L_\infty( {\cal E}_p) = \max_{k=1,2,3}  \,\,\,     \max_{i  \in \{ I \}    }   | {\bf w}_i(k) - {\bf w}^{exact}_i(k) |.
\end{eqnarray}

\begin{figure}[th!]
    \centering
      \hfill  
                \begin{subfigure}[t]{0.48\textwidth}
        \includegraphics[width=0.9\textwidth]{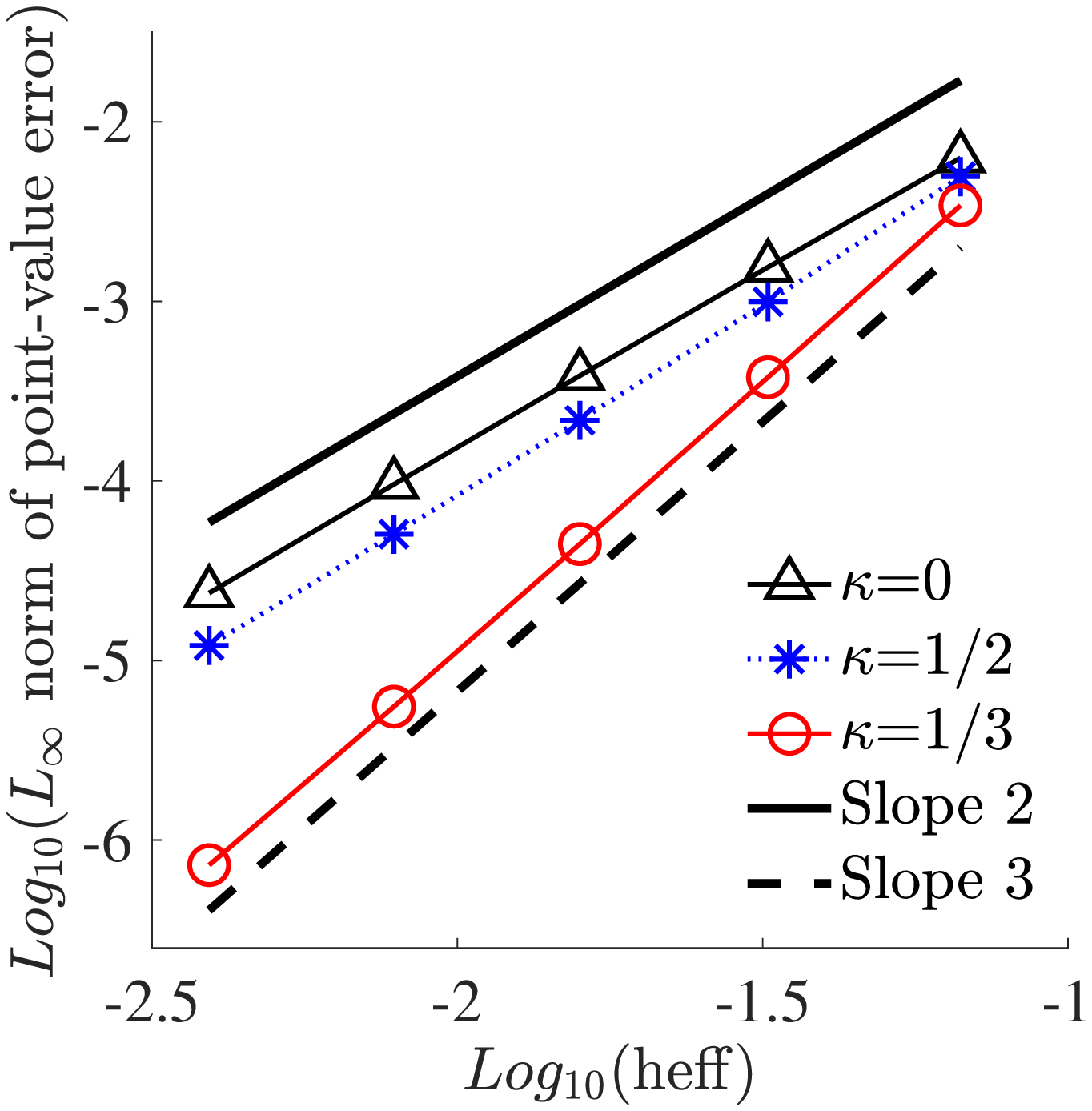}
          \caption{\color{black}  $c_\epsilon = 0.05$.}
       \label{fig:euler_1d_case1}
      \end{subfigure}
      \hfill
          \begin{subfigure}[t]{0.48\textwidth}
        \includegraphics[width=0.9\textwidth]{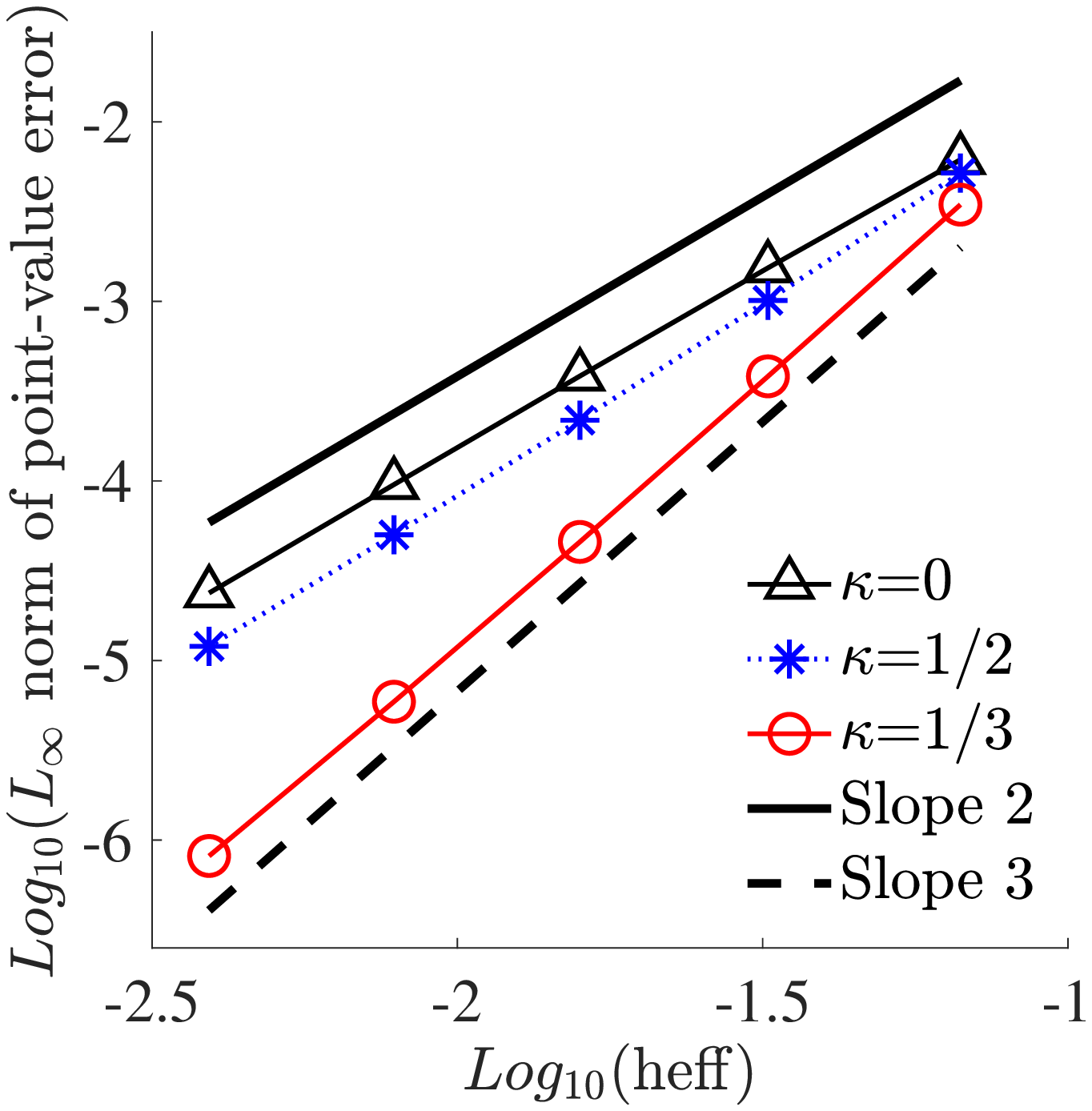}
          \caption{ \color{black}  $c_\epsilon = 0.1$.}
       \label{fig:euler_1d_case2}
      \end{subfigure}
      \hfill
          \begin{subfigure}[t]{0.48\textwidth}
        \includegraphics[width=0.9\textwidth]{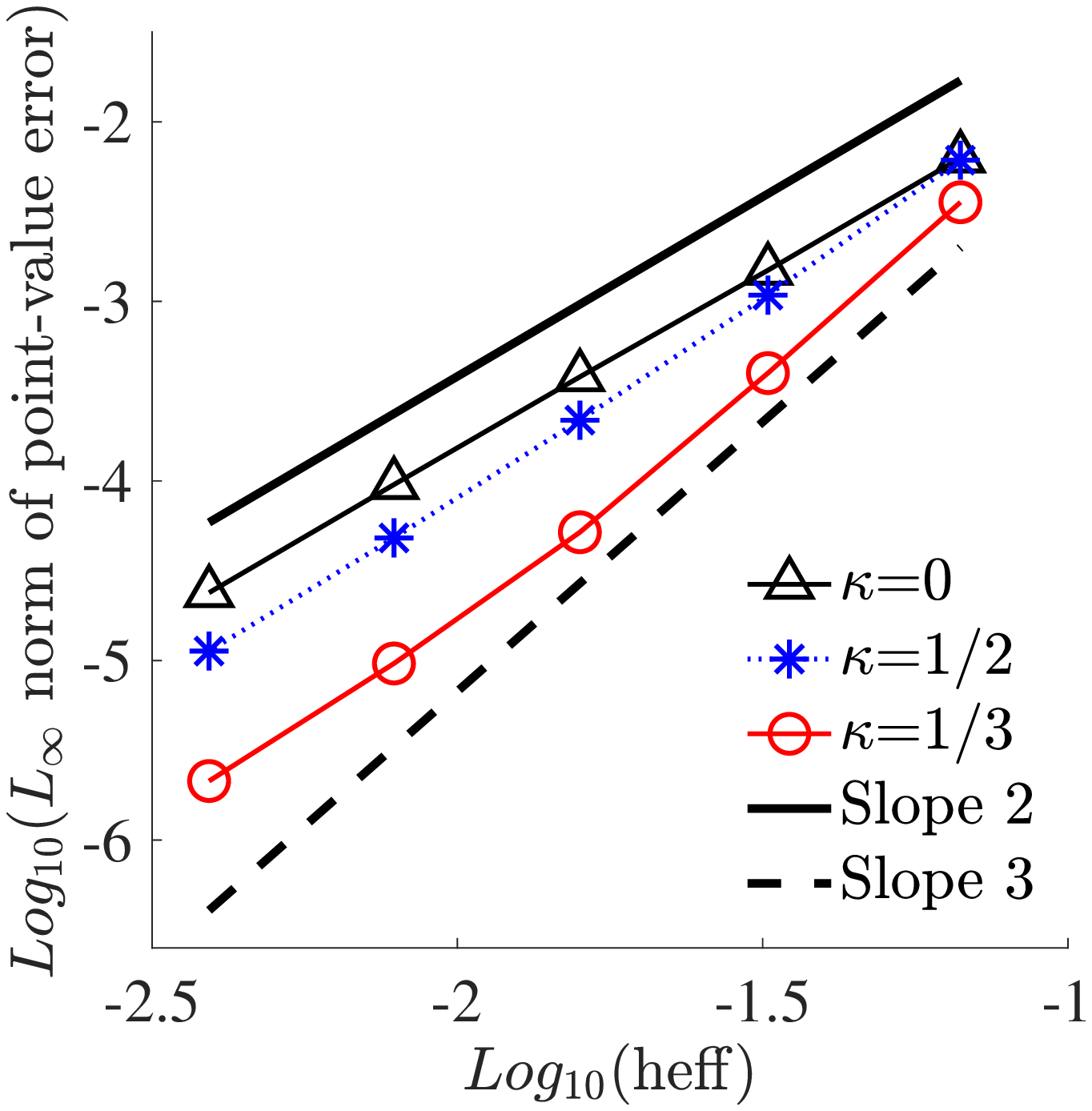}
          \caption{\color{black}   $c_\epsilon = 0.3$.}
       \label{fig:euler_1d_case3}
      \end{subfigure}
      \hfill  
                \begin{subfigure}[t]{0.48\textwidth}
        \includegraphics[width=0.9\textwidth]{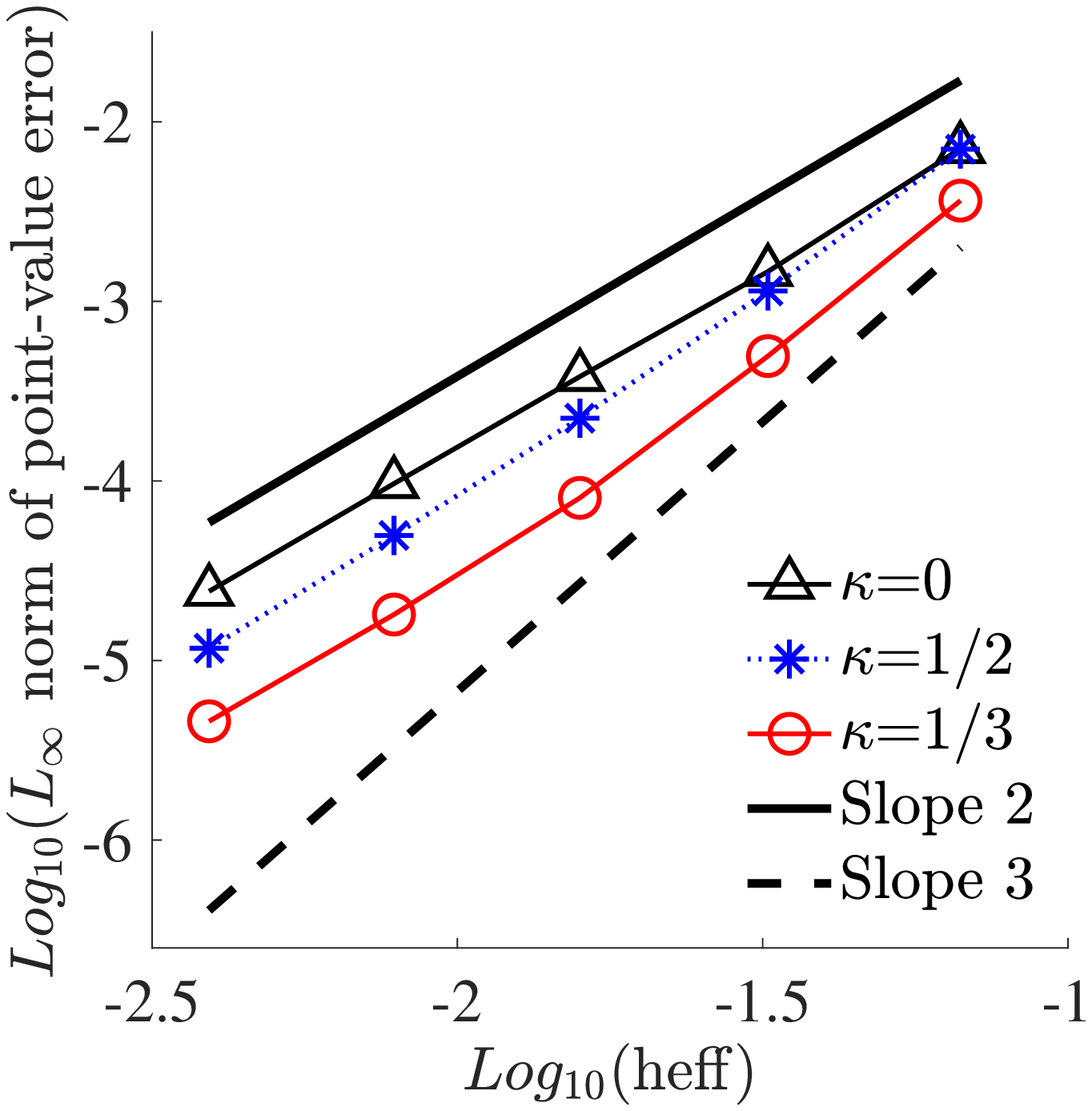}
          \caption{\color{black}   $c_\epsilon = 0.5$.}
       \label{fig:euler_1d_case4}
      \end{subfigure}
      \hfill    
      \\
            \caption{
\label{fig:euler_1d_case}%
UMUSCL: error convergence for the Euler equations in one dimension.
.
} 
\end{figure}

Results are shown in Figure \ref{fig:euler_1d_case}. {\color{black} As expected, third-order accuracy is observed for $c_\epsilon=0.05$ and $0.1$, when $\kappa=1/3$ as shown in Figures \ref{fig:euler_1d_case1} and \ref{fig:euler_1d_case2}. For $c_\epsilon=0.3$ and $0.5$, the second-order error begins to dominate fairly accurately as predicted by the estimates $h_{cr}$, as can be seen in Figures \ref{fig:euler_1d_case3} and \ref{fig:euler_1d_case4}. 
From the fact that the actual value of $\epsilon$ is $0.15$ for the last case ($c_\epsilon=0.5$), we would expect that 
false third-order accuracy can be entirely avoided if $\epsilon$ is set equal to the corresponding coefficients for the other variables, i.e., $\epsilon = 0.2$. 
 }
%


   

\begin{figure}[th!]
    \centering
      \hfill  
                \begin{subfigure}[t]{0.48\textwidth}
        \includegraphics[width=0.9\textwidth]{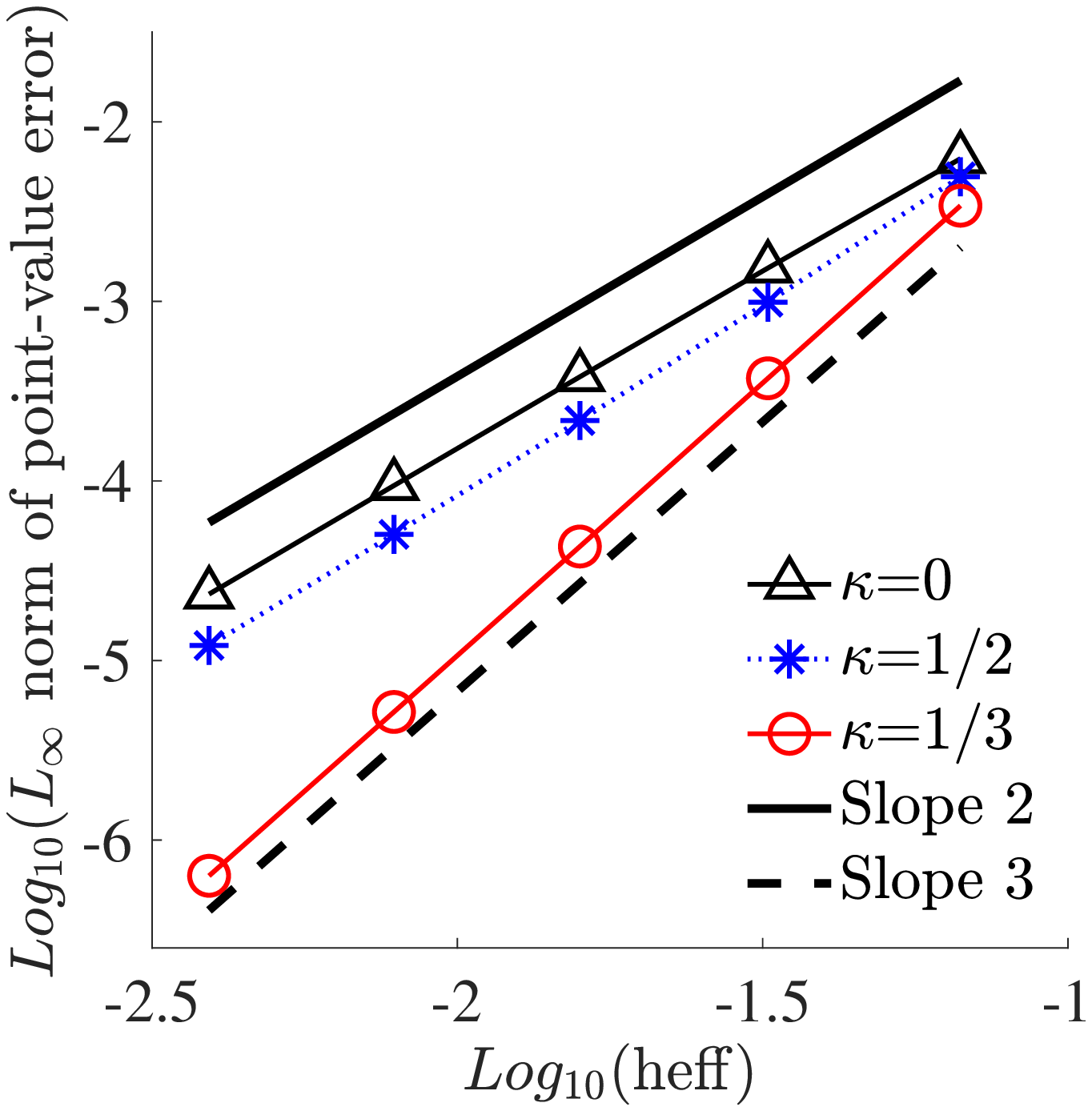}
          \caption{\color{black}  $c_\epsilon = 0.05$.}
       \label{fig:euler_1d_case1}
      \end{subfigure}
      \hfill
          \begin{subfigure}[t]{0.48\textwidth}
        \includegraphics[width=0.9\textwidth]{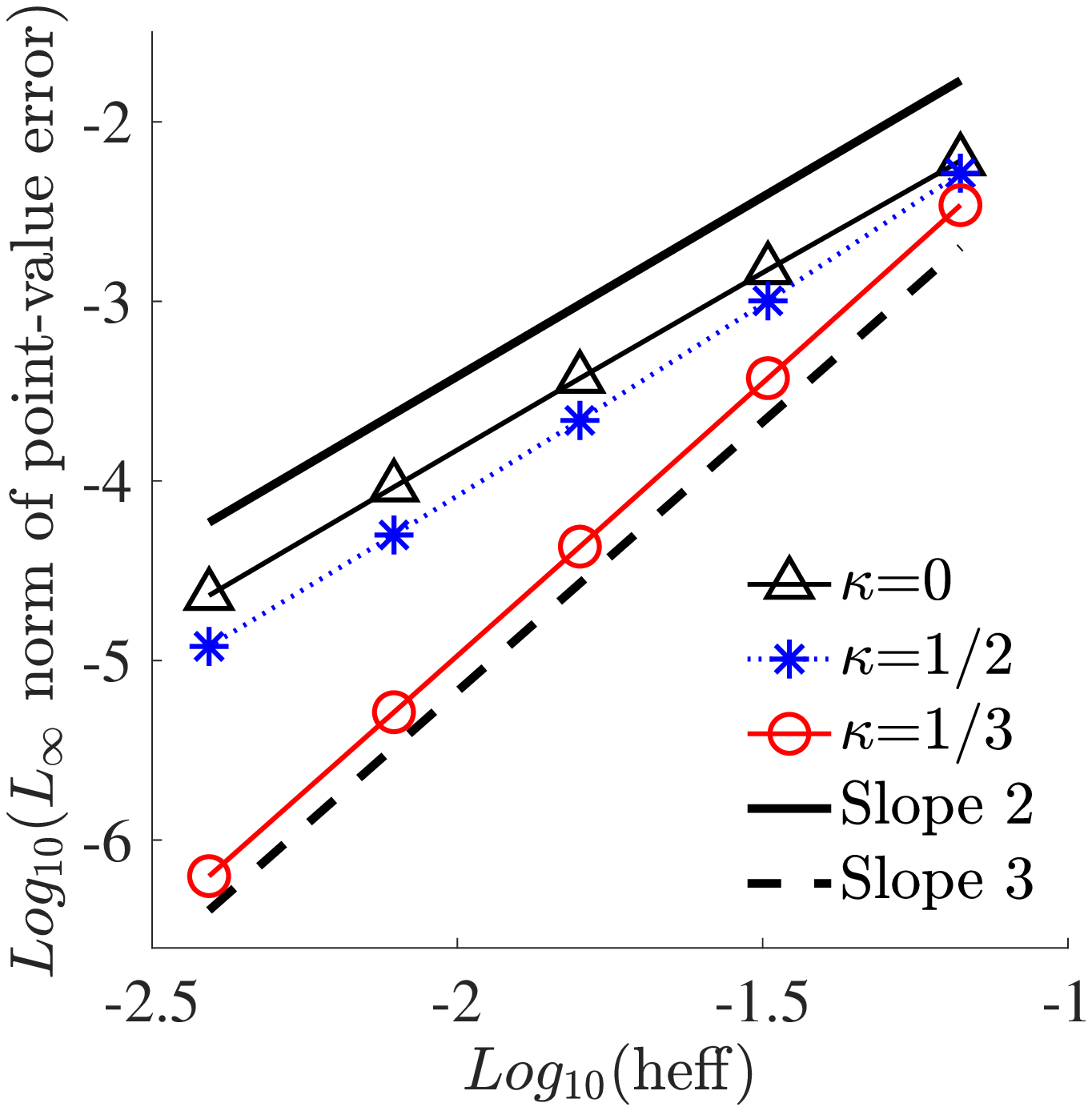}
          \caption{\color{black}  $c_\epsilon = 0.1$.}
       \label{fig:euler_1d_case2}
      \end{subfigure}
      \hfill
          \begin{subfigure}[t]{0.48\textwidth}
        \includegraphics[width=0.9\textwidth]{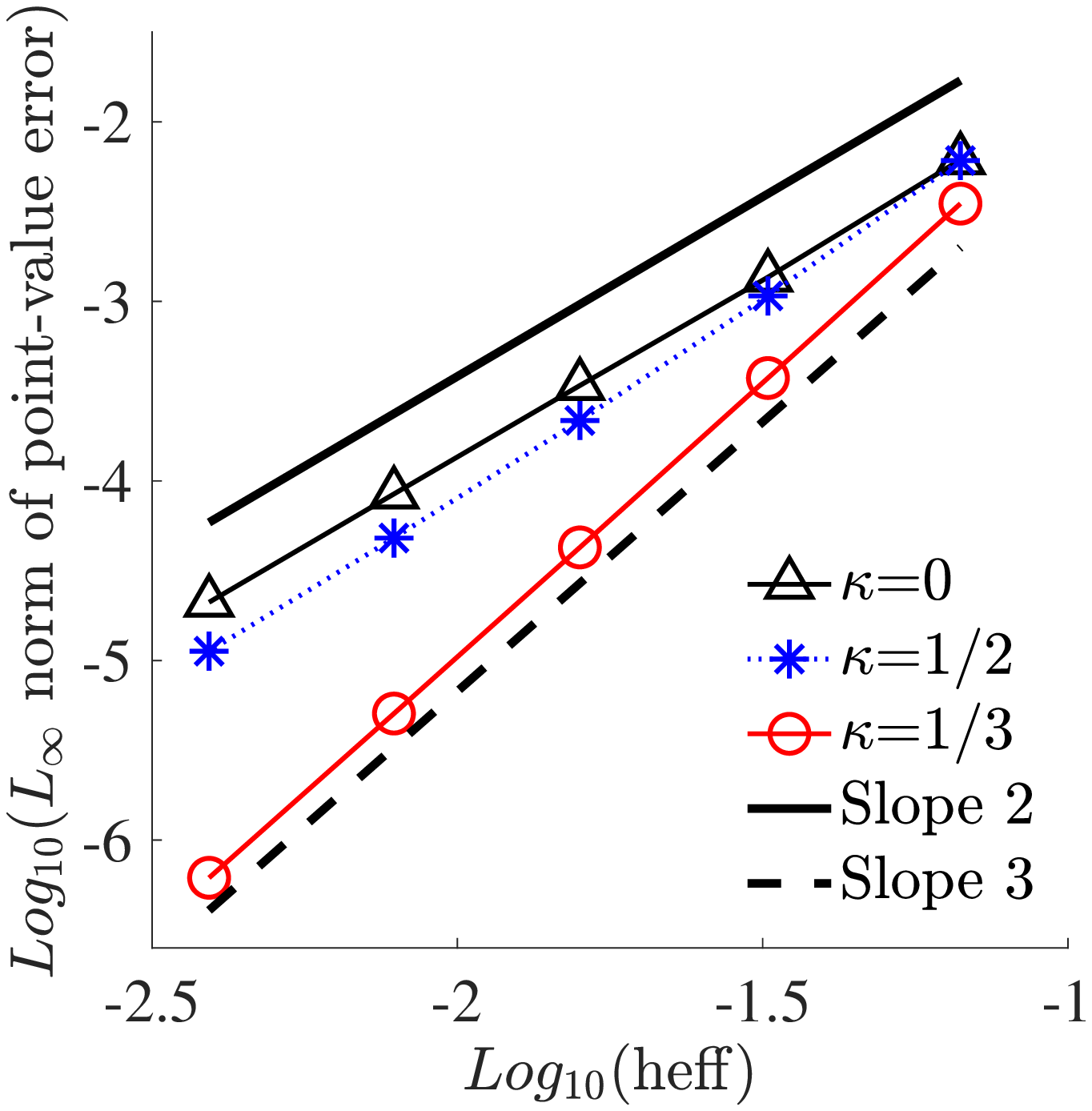}
          \caption{\color{black}  $c_\epsilon = 0.3$.}
       \label{fig:euler_1d_case3}
      \end{subfigure}
      \hfill  
                \begin{subfigure}[t]{0.48\textwidth}
        \includegraphics[width=0.9\textwidth]{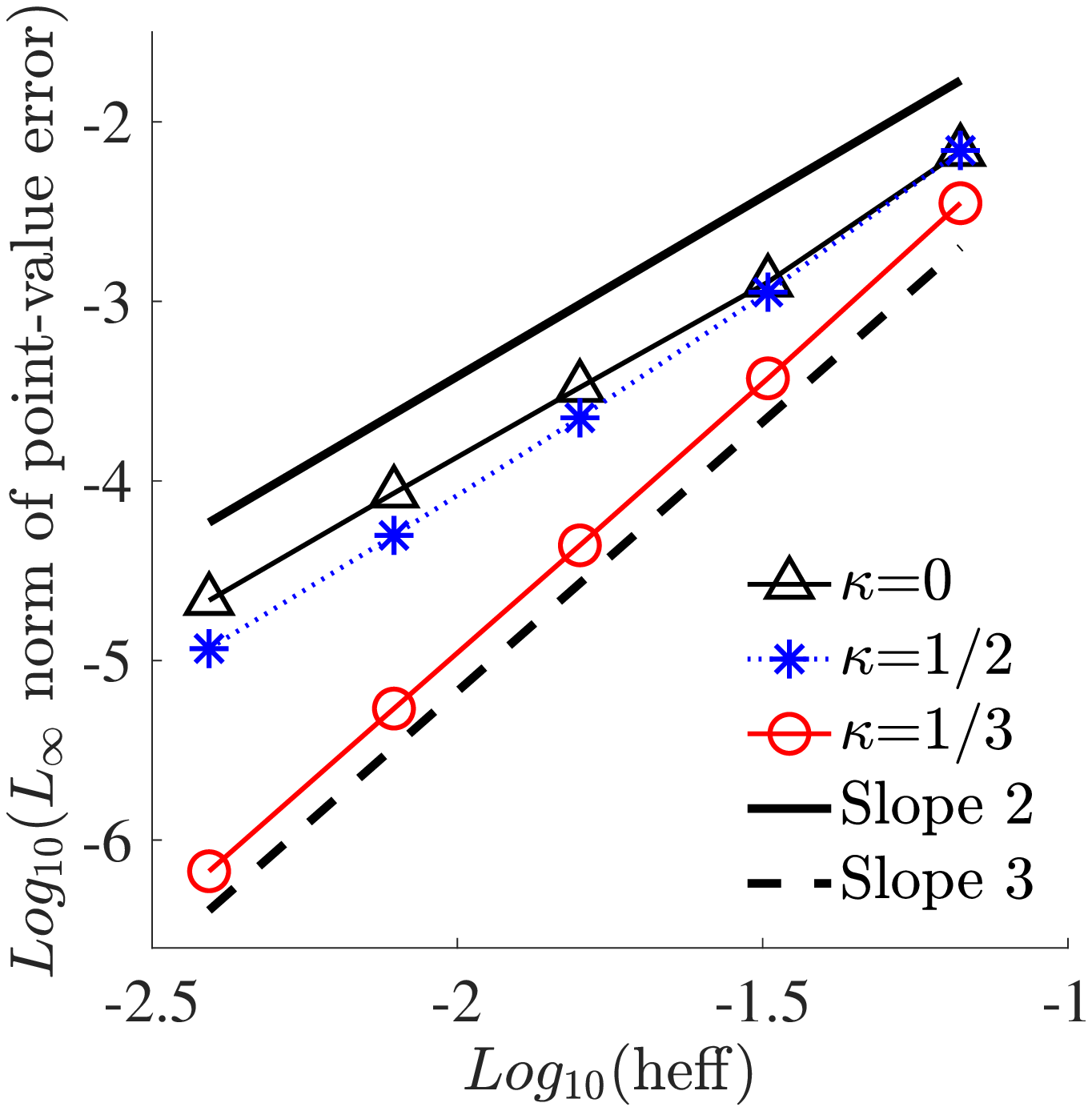}
          \caption{ \color{black} $c_\epsilon = 0.5$.}
       \label{fig:euler_1d_case4}
      \end{subfigure}
      \hfill    
      \\
            \caption{
\label{fig:euler_FR_1d_case}%
FSR: error convergence for the Euler equations in one dimension.
} 
\end{figure}

To demonstrate that the UMUSCL scheme can be made genuinely third-order with the flux reconstruction (i.e., FSR), we consider the following numerical flux,
\begin{eqnarray}
{\bf F}({\bf w}_L,{\bf w}_R, {\bf f}_L , {\bf f}_R) =  \frac{1}{2} \left[   {\bf f}_L   +  {\bf f}_R  \right]  - \frac{\hat{\bf D}}{2}  \left[  {\bf u}({\bf w}_R) -  {\bf u}({\bf w}_L) \right],
\label{cfd_numerical_flux_euler}
\end{eqnarray}
where the left and right fluxes ${\bf f}_L $ and ${\bf f}_R$ are directly reconstructed as in Equations (\ref{flux_reconstruct_fL}) and (\ref{flux_reconstruct_fR}). 
This is slightly more expensive because ${\bf w}_L$ and ${\bf w}_R$ still need to be computed as well, for example. However, the scheme is genuinely third-order 
and thus third-order accuracy is obtained for any value of $c_\epsilon$ as can be seen in Figure \ref{fig:euler_FR_1d_case}. 


\begin{figure}[th!]
    \centering
      \hfill  
                \begin{subfigure}[t]{0.32\textwidth}
        \includegraphics[width=0.9\textwidth]{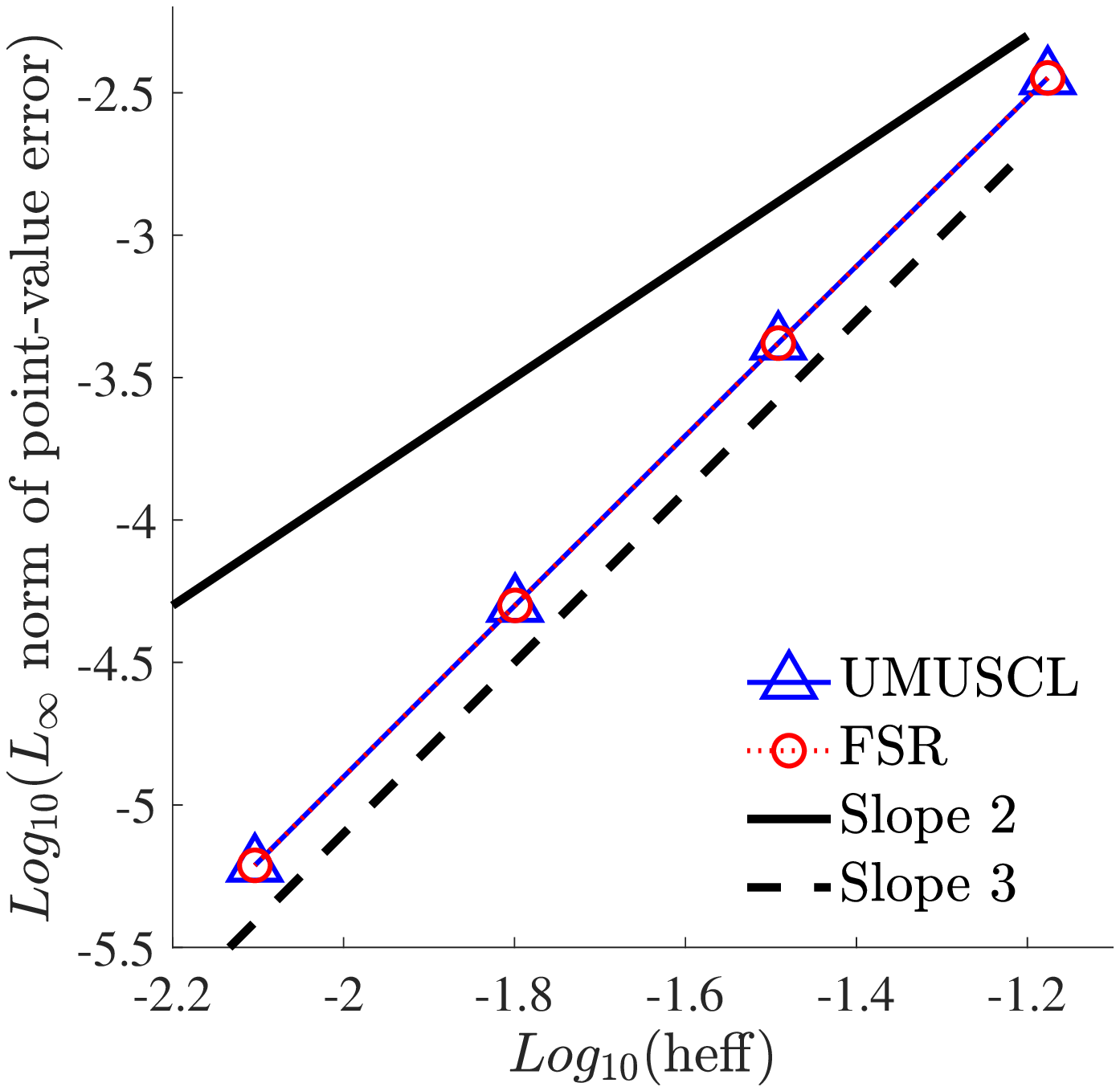}
          \caption{$\rho=\rho(x)$, $u=0.3$, $p=1$}
       \label{fig:euler_1d_case100}
      \end{subfigure}
      \hfill
          \begin{subfigure}[t]{0.32\textwidth}
        \includegraphics[width=0.9\textwidth]{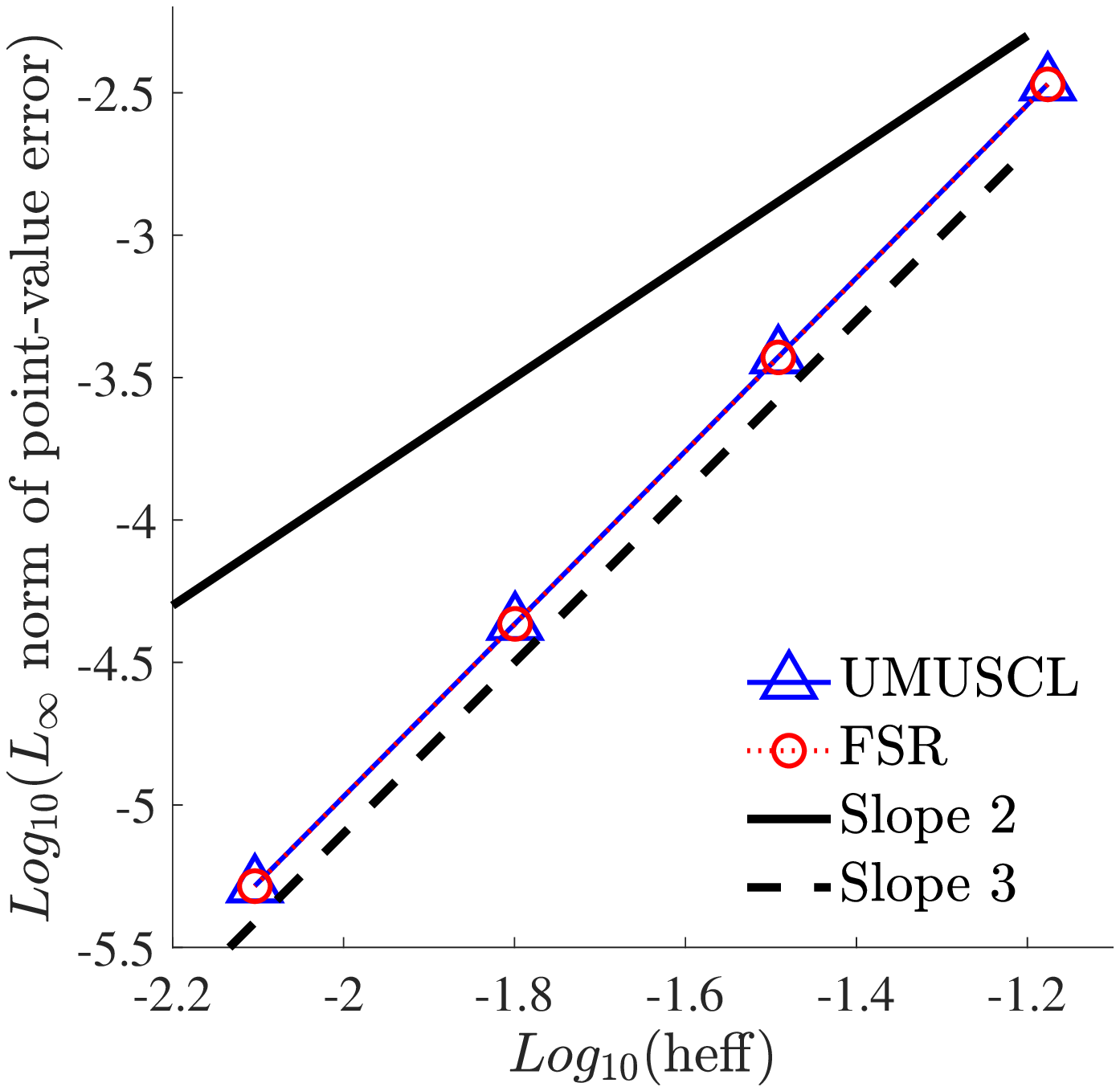}
          \caption{$\rho=\rho(x)$, $u=0.3$, $p=p(x)$}
       \label{fig:euler_1d_case101}
      \end{subfigure}
      \hfill
          \begin{subfigure}[t]{0.32\textwidth}
        \includegraphics[width=0.9\textwidth]{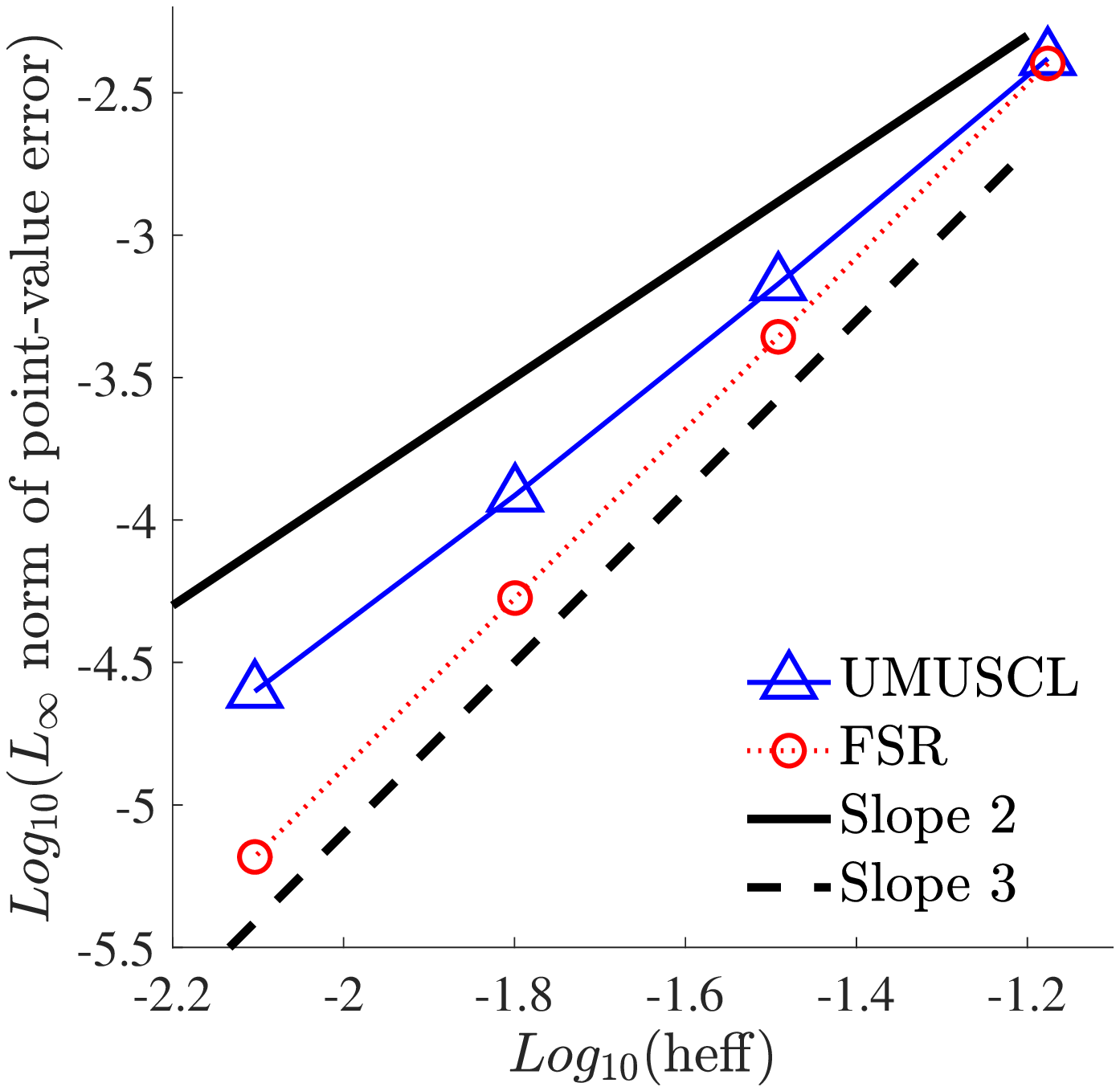}
          \caption{$\rho=\rho(x)$, $u=u(x)$, $p=p(x)$}
       \label{fig:euler_1d_case102}
      \end{subfigure}
      \hfill    
      \\
            \caption{
\label{fig:euler_FR_1d_case_simple}%
UMUSCL and FSR: error convergence for the Euler equations in one dimension: 
(a) $\epsilon_\rho=0.2$, $\epsilon=\epsilon_p=0$, 
(b) $\epsilon_\rho=\epsilon_p=0.2$, $\epsilon=0$, 
(c) $\epsilon_\rho=\epsilon_p=\epsilon=0.2$
.
} 
\end{figure}

For the sake of completeness, we present results for the following three cases: 
(a)  $\epsilon_\rho=0.2$ and $\epsilon=\epsilon_p=0$ (an entropy wave solution \cite{ZhongSheng:CF2020}),  
(b)  $\epsilon_\rho=\epsilon_p=0.2$ and $\epsilon=0$,
(c)  $\epsilon_\rho=\epsilon_p=\epsilon=0.2$. 
As mentioned in the previous section, the Euler equations are linearized in the cases (a) and (b), where the velocity is constant. Therefore,
the UMUSCL scheme may give third-order accuracy. However, the case (c) renders the Euler equations fully nonlinear and thus 
is expected to reveal second-order accuracy of the UMUSCL scheme. For all the cases, UMUSCL is compared with FSR and $\kappa=1/3$ is used for both schemes.
Results are shown in Figure \ref{fig:euler_FR_1d_case_simple}. As expected, UMUSCL gives third-order accuracy for 
 the cases (a) and (b), but deteriorates to second-order in the case (c). 
On the other hand, FSR gives third-order accuracy for all the cases. 

\subsection{Euler equations in two dimensions}
\label{results_euler_2d} 

Finally, we consider two-dimensional cases for the Euler equations, which have been used to mistakenly verify high-order accuracy of the UMUSCL scheme 
in Refs.\cite{yang_harris:AIAAJ2016,yang_harris:CCP2018,Nishikawa_FANG_AQ:Aviation2020} (and also of a relevant scheme in Ref.\cite{ZhongSheng:CF2020}). 
We will demonstrate that third-order accuracy is due to the unexpected linearization and thus not genuine but that 
genuine third-order accuracy can be achieved with the flux reconstruction.

Consider the Euler equations in two dimensions, 
\begin{eqnarray}
\frac{ \partial {\bf u}}{\partial t} + \mbox{div} \, {\cal F} 
  = {\bf s} ,
  \quad
  {\bf u} =  \left[  \begin{array}{c} 
               \rho       \\ [1ex]
               \rho {\bf v}     \\ [1ex]
               \rho E     
              \end{array} \right], \quad
  {\cal F} =  
 \left[  \begin{array}{c} 
               \rho {\bf v}^t      \\  [1ex]
               \rho  {\bf v} \dyad  {\bf v}  + p  {\bf I }  \\ [1ex]
               \rho  {\bf v}^t H 
              \end{array} \right] ,
              \label{2D_NS_Conservative}
\end{eqnarray}
where  ${\bf v}$ is the velocity vector with Cartesian components ${\bf v} = (u,v)^t$ (which is a column vector; the superscript indicates the transpose), $ \dyad$ denotes the dyadic product, ${\bf I}$ is the identity matrix, 
and $E = H-p / \rho  = (p/\rho) / (\gamma-1)  + {\bf v}^2 /2$ is the specific total energy.

  \begin{figure}[t]
    \centering 
        \includegraphics[width=0.48\textwidth]{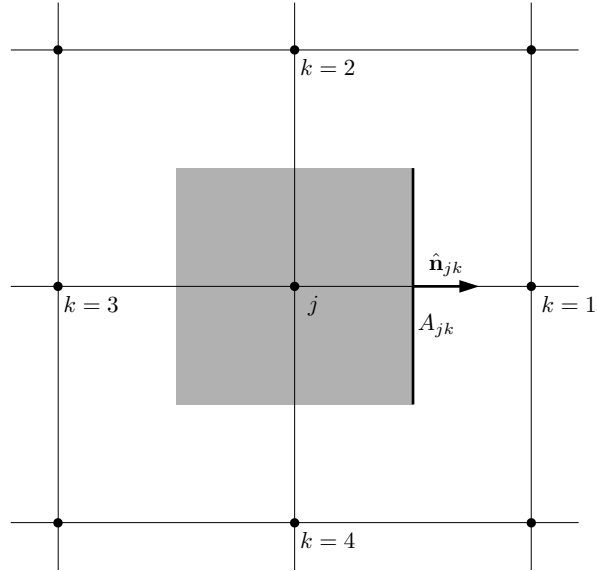} 
            \caption{
\label{fig:oned_fv_data_quad_inter_uns}%
A control volume around a node $j$ with local numbering of the neighbor nodes typical in unstructured grids.
} 
\end{figure}

The UMUSCL scheme is implemented in the form applicable to a general unstructured grid \cite{yang_harris:AIAAJ2016}, i.e., at a node $j$ as in 
Figure \ref{fig:oned_fv_data_quad_inter_uns}, 
\begin{eqnarray}
\frac{d {\bf u}_{j} }{dt}   + \frac{1}{V_j}   \sum_{ k \in \{ k_j \}} {\boldsymbol \Phi}_{jk}(u_L,u_R)  A_{jk}= {\bf s}_{j},
\label{semi_discrete_general_2d_euler}
\end{eqnarray}
with the Roe flux,
\begin{eqnarray}
{\bf F}({\bf w}_L,{\bf w}_R) =  \frac{1}{2} \left[   {\bf f}({\bf w}_L)   +  {\bf f}({\bf w}_R)   \right]  - \frac{\hat{\bf D}_n}{2}  \left[  {\bf u}({\bf w}_R) -  {\bf u}({\bf w}_L) \right],
\label{cfd_numerical_flux_2d_euler}
\end{eqnarray}
where $ {\bf f} = {\cal F} \cdot \hat{\bf n}_{jk}$ is the flux projected along the face normal ${\bf n}_{jk} = \hat{\bf n}_{jk} |{\bf n}_{jk}| = \hat{\bf n}_{jk}  A_{jk} $, $\hat{\bf D}_n = | \partial {\bf f} / \partial {\bf u}  |$ is the dissipation term evaluated with the Roe averages, the solution reconstruction is performed with the primitive variables 
${\bf w} = (\rho,u,v,p)$:
\begin{eqnarray}
{\bf w}_L &=& \kappa  \frac{ {\bf w}_{j} +  {\bf w}_{k}}{2} +   (1-\kappa) \left[  {\bf w}_j  + \frac{1}{2}  \nabla {\bf w}_j \cdot (  {\bf x}_k - {\bf x}_j  )  \right] ,  
\label{umuscl_L2_2d_euler}\\ [1ex]
{\bf w}_R &=& \kappa  \frac{ {\bf w}_{k} + {\bf w}_{j}}{2} +   (1-\kappa) \left[  {\bf w}_{k}  +  \frac{1}{2}  \nabla {\bf w}_k \cdot (  {\bf x}_j - {\bf x}_k  )   \right], 
 \label{umuscl_R2_2d_euler}
\end{eqnarray}
where ${\bf x}_j$ and ${\bf x}_k$ denote the nodal coordinates of $j$ and its neighbor $k$, respectively and the gradients are computed by a linear least-squares
method. On Cartesian grids, the scheme (\ref{semi_discrete_general_2d_euler}) reduces to the form (\ref{semi_discrete}) with the $\kappa$-reconstruction scheme applied as a one-dimensional algorithm in each coordinate direction. Note that no high-order flux quadrature is used; the flux is computed only at the midpoint of each face. 

For the flux reconstruction version, FSR, the numerical flux is computed as
\begin{eqnarray}
{\bf F}({\bf w}_L,{\bf w}_R, {\bf f}_L , {\bf f}_R) =  \frac{1}{2} \left[   {\bf f}_L  +  {\bf f}_R   \right]  - \frac{\hat{\bf D}_n}{2}  \left[  {\bf u}({\bf w}_R) -  {\bf u}({\bf w}_L) \right],
\label{cfd_numerical_flux_euler_2d_euler_FR}
\end{eqnarray}
where
\begin{eqnarray}
{\bf f}_L &=& \kappa  \frac{ {\bf f}({\bf w}_j) + {\bf f}({\bf w}_k)   }{2} +   (1-\kappa) \left[ {\bf f}({\bf w}_i)  + \frac{1}{2}  \left(  \frac{\partial {\bf f}}{\partial {\bf w}} \right)_{\!\! j}  \nabla {\bf w}_j \cdot (  {\bf x}_k - {\bf x}_j  )  \right] ,  
\label{umuscl_L2_2d_euler_FR}\\ [1ex]
{\bf f}_R &=& \kappa  \frac{ {\bf f}({\bf w}_k) +{\bf f}({\bf w}_j) }{2} +   (1-\kappa) \left[  {\bf f}({\bf w}_k) +  \frac{1}{2}   \left(  \frac{\partial {\bf f}}{\partial {\bf w}} \right)_{\!\! k}  \nabla {\bf w}_k \cdot (  {\bf x}_j - {\bf x}_k  )   \right], 
 \label{umuscl_R2_2d_euler_FR}
\end{eqnarray}
and the flux Jacobian is given, with the unit face normal vector $\hat{\bf n}_{jk}$ and ${\bf v}$ defined as column vectors and the notation $u_n = {\bf v} \cdot  {\bf n}_{jk}$, as
\begin{eqnarray}
 \frac{\partial {\bf f}}{\partial {\bf w}} 
 =
  \left[  \begin{array}{ccc} 
u_n          &  \rho  \hat{\bf n}_{jk}^t   & 0    \\ [1ex]
u_n {\bf v} & \rho ( u_n {\bf I} + {\bf v} \dyad   \hat{\bf n}_{jk} )  &  \hat{\bf n}_{jk}    \\ [1ex]
u_n {\bf v}^2 /2&  \rho(H  \hat{\bf n}_{jk}^t + u_n {\bf v}^t ) & \gamma u_n /(\gamma-1)
              \end{array} \right].
\end{eqnarray}
Note that the flux gradients are efficiently computed from the solution gradients by the chain rule. Therefore, it is not necessary to compute the flux gradients directly by a least-squares method and thus does not also require extra storage for flux gradients. This efficient version of FSR is referred to as FSR-CR (chain rule) in the rest of the paper. 
Ref.\cite{Barakos:IJNMF2018} seems to perform the flux reconstruction and therefore their method may be genuinely high-order accurate. However, their method requires the computation and storage of the the flux gradients; it can be very expensive in three dimensions. A further discussion including fourth- and fifth-order accurate schemes will be given in a subsequent paper.

 


\subsubsection{Steady problem}
\label{results_euler_2d_steady} 

\begin{figure}[htbp!]
    \centering
      \begin{subfigure}[t]{0.32\textwidth}
        \includegraphics[width=\textwidth,trim=0 0 0 0,clip]{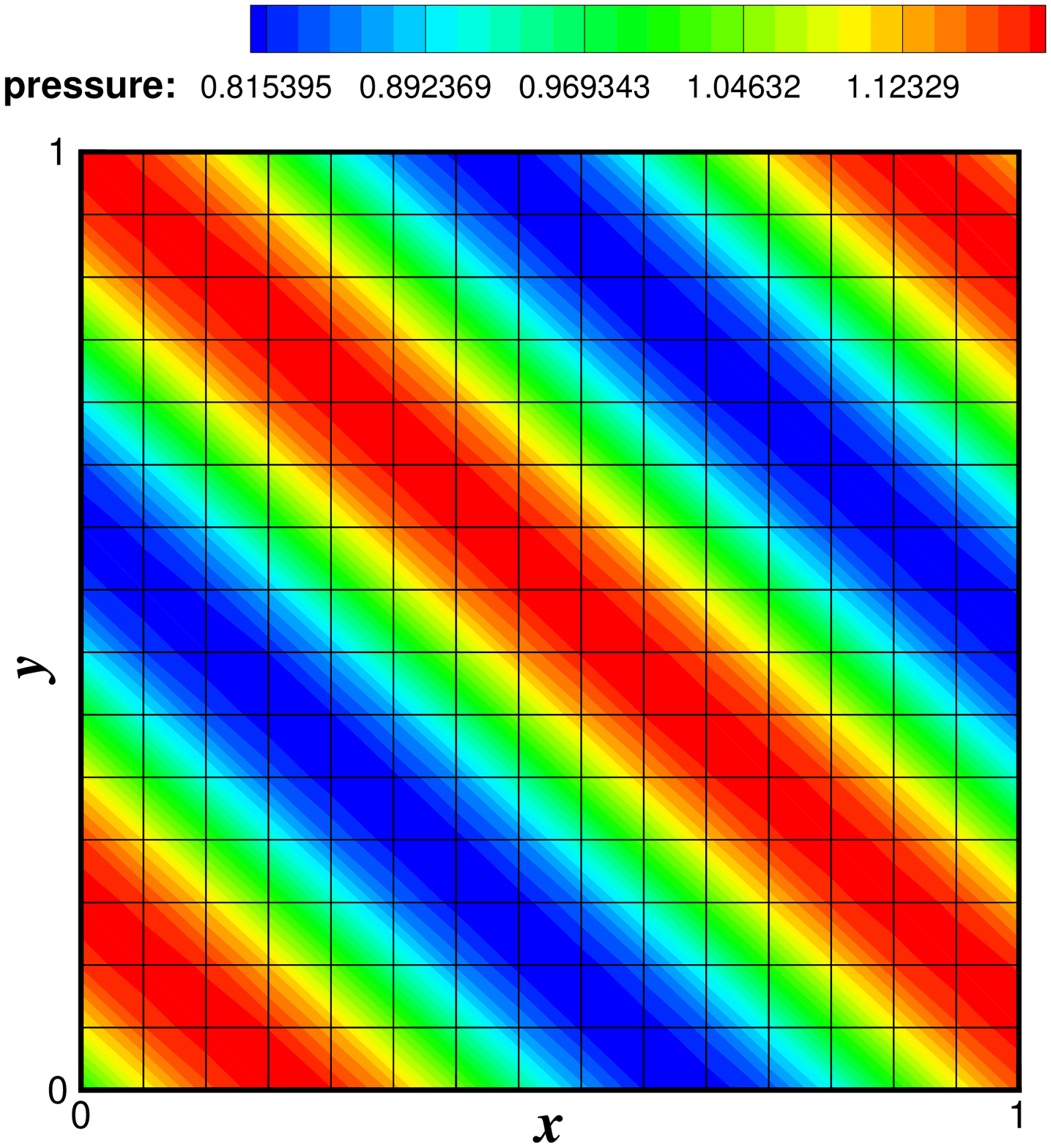}
          \caption{The coarsest grid and solution.}
          \label{fig:inv_err_mms_gridsol}
      \end{subfigure}
      \hfill
      \begin{subfigure}[t]{0.32\textwidth}
        \includegraphics[width=\textwidth,trim=0 0 0 0,clip]{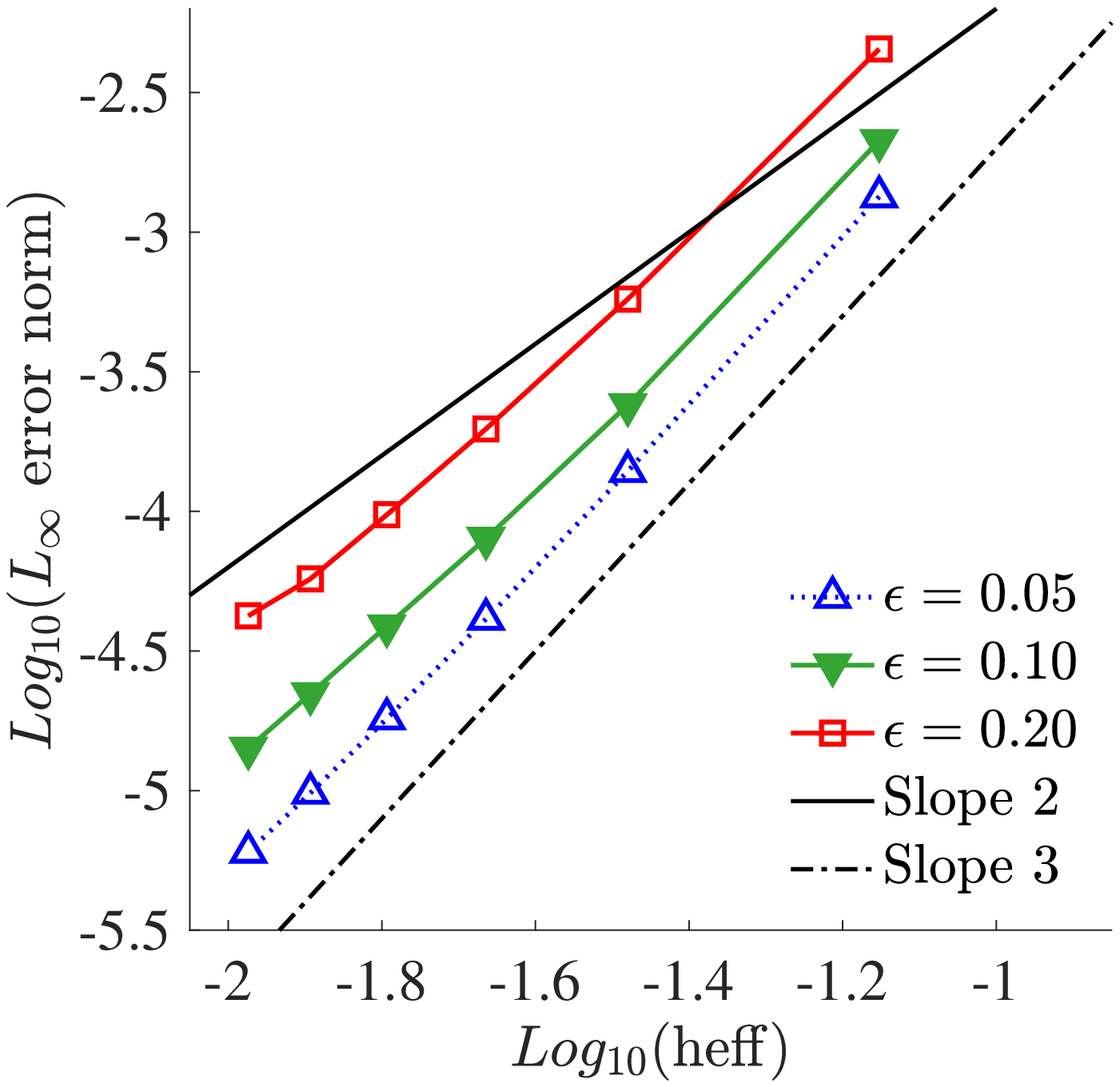}
          \caption{UMUSCL}
          \label{fig:inv_err_mms_error}
      \end{subfigure}
      \hfill
      \begin{subfigure}[t]{0.32\textwidth}
        \includegraphics[width=\textwidth,trim=0 0 0 0,clip]{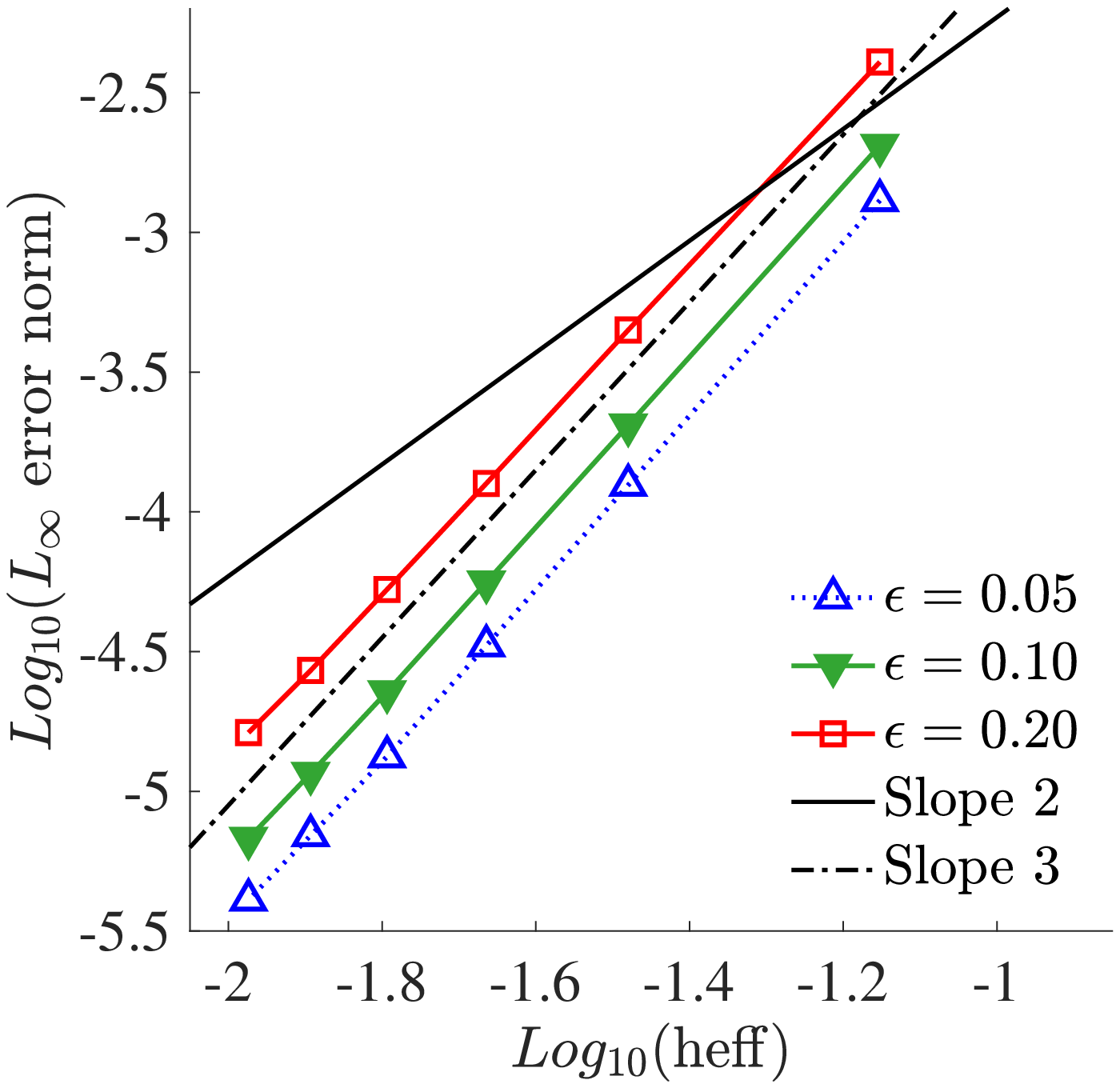}
          \caption{FSR-CR}
          \label{fig:inv_err_mms_error_fr}
      \end{subfigure}
      \caption{Error convergence study for the steady Euler equations in two dimensions.}
\label{fig:inv_err_mms} 
\end{figure}

Let us begin with a steady case in a unit square, similar to the one used in Ref.\cite{Nishikawa_FANG_AQ:Aviation2020}, with the exact solutions defined by
\begin{eqnarray}
\rho                &=&  1  + 0.2  \sin[   \pi (  2.3  x + 2.3 y ) ], \\
u                    &=&  u_\infty  + \epsilon \sin[    \pi ( 2  x + 2 y )  ], \\ 
v                    &=&  v_\infty  + \epsilon \sin[    \pi ( 2   x + 2 y ) ], \\ 
p                   &=&  1  + 0.2 \sin[     \pi ( 2.5  x + 2.5 y ) ],
\end{eqnarray}
where $u_\infty  = 0.15 $, $v_\infty = 0.02$, and $\epsilon$ is given the following values:
\begin{eqnarray}
 \epsilon   = 0, \quad 0.05, \quad 0.1, \quad 0.2.
\end{eqnarray}
The forcing term ${\bf s}_j$ is numerically computed with the above exact solutions at a node $j$ and added to the residual (see Ref.\cite{nishikawa_centroid:JCP2020} for details). 
To exclude boundary effects, we specify the exact solution at nodes on the boundary and at their neighbors. 
The residual equations are solved using an 
implicit defect-correction solver (see, e.g., Ref.\cite{NakashimaWatanabeNishikawa_AIAA2016-1101}) until the residual norm is reduced
 by seven orders of magnitude for a series of regular quadrilateral grids: 49$\times$49, 65$\times$65, 81$\times$81,
 97$\times$97, 113$\times$113, 129$\times$129. Figure \ref{fig:inv_err_mms_gridsol} shows exact pressure contours on a
 coarser 17$\times$17 grid for the sake of illustration. For this problem, we compare UMUSCL and FSR-CR with $\kappa=1/3$.

%

Figure \ref{fig:inv_err_mms_error} shows the $L_\infty$ error convergence results for the pressure obtained with the UMUSCL scheme.
As expected, the convergence rate deteriorates as $\epsilon$ increases towards $0.2$. 
Nearly third-order convergence is observed 
for $\epsilon=0.05$, but it eventually deteriorates to second-order at $\epsilon=0.2$. 
 The exact solution used in Ref.\cite{Nishikawa_FANG_AQ:Aviation2020} is similar to the above solution with $\epsilon=0.05$; thus third-order accuracy reported in 
 Ref.\cite{Nishikawa_FANG_AQ:Aviation2020} is not genuine (the author did not realize it at the time of submitting the paper). It is important to observe that the convergence rates are all close to third-order on coarse grids in Figure \ref{fig:inv_err_mms_error}: accuracy verification must be conducted with sufficiently fine grids to observe the true order of accuracy. 
  On the other hand, if we perform the flux reconstruction, then third-order accuracy is achieved for all values of $\epsilon$ as 
 shown in Figure \ref{fig:inv_err_mms_error_fr}. 

\subsubsection{Unsteady inviscid vortex transport}
\label{results_euler_2d_unsteady} 

\begin{figure}[htbp!]
    \centering
      \begin{subfigure}[t]{0.32\textwidth}
        \includegraphics[width=\textwidth,trim=0 0 0 0,clip]{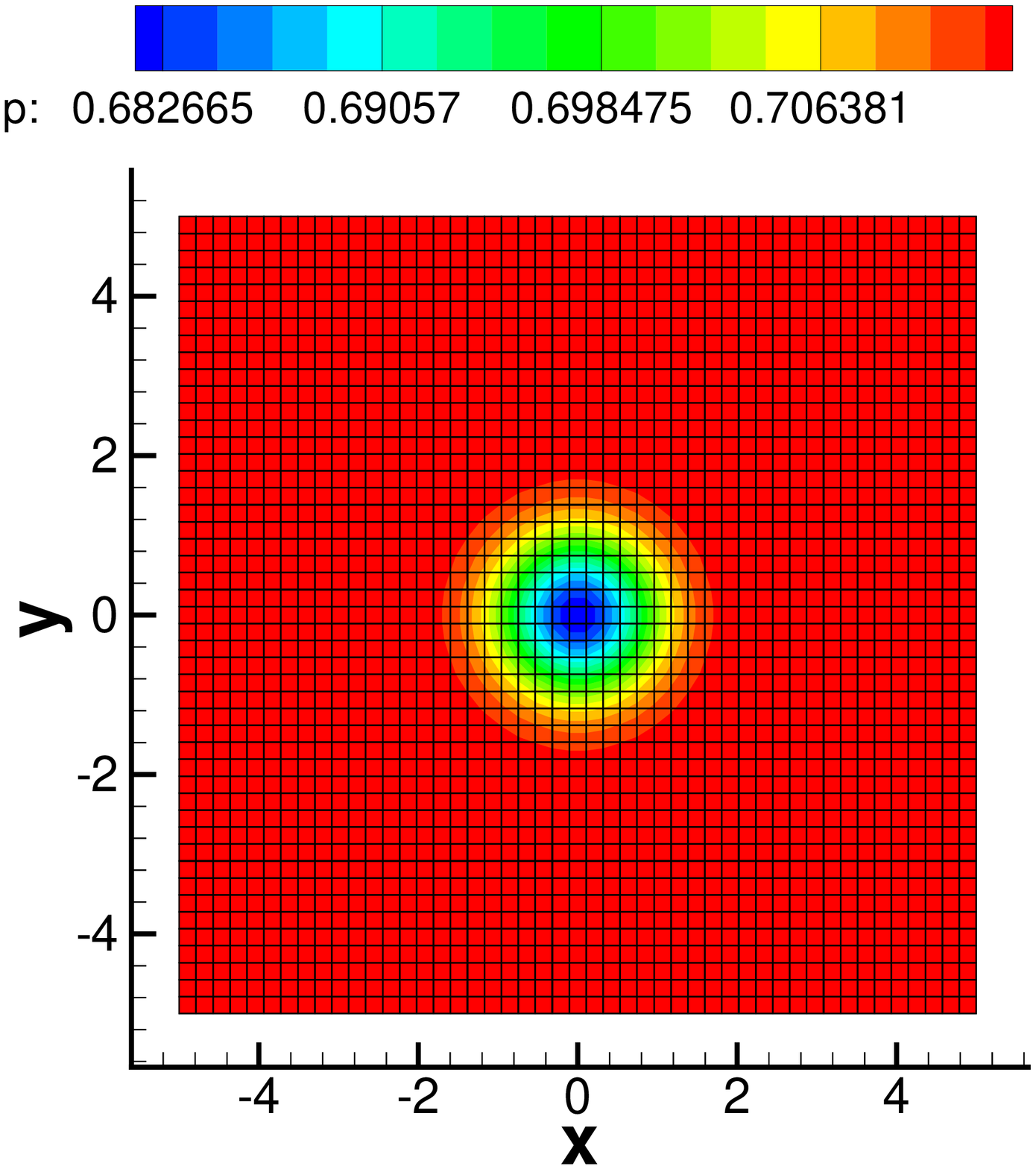}
          \caption{Initial solution on the coarsest grid ($K=1$).}
          \label{fig:inv_uns_err_mms_gridsol}
      \end{subfigure}
      \hfill
      \begin{subfigure}[t]{0.32\textwidth}
        \includegraphics[width=\textwidth,trim=0 0 0 0,clip]{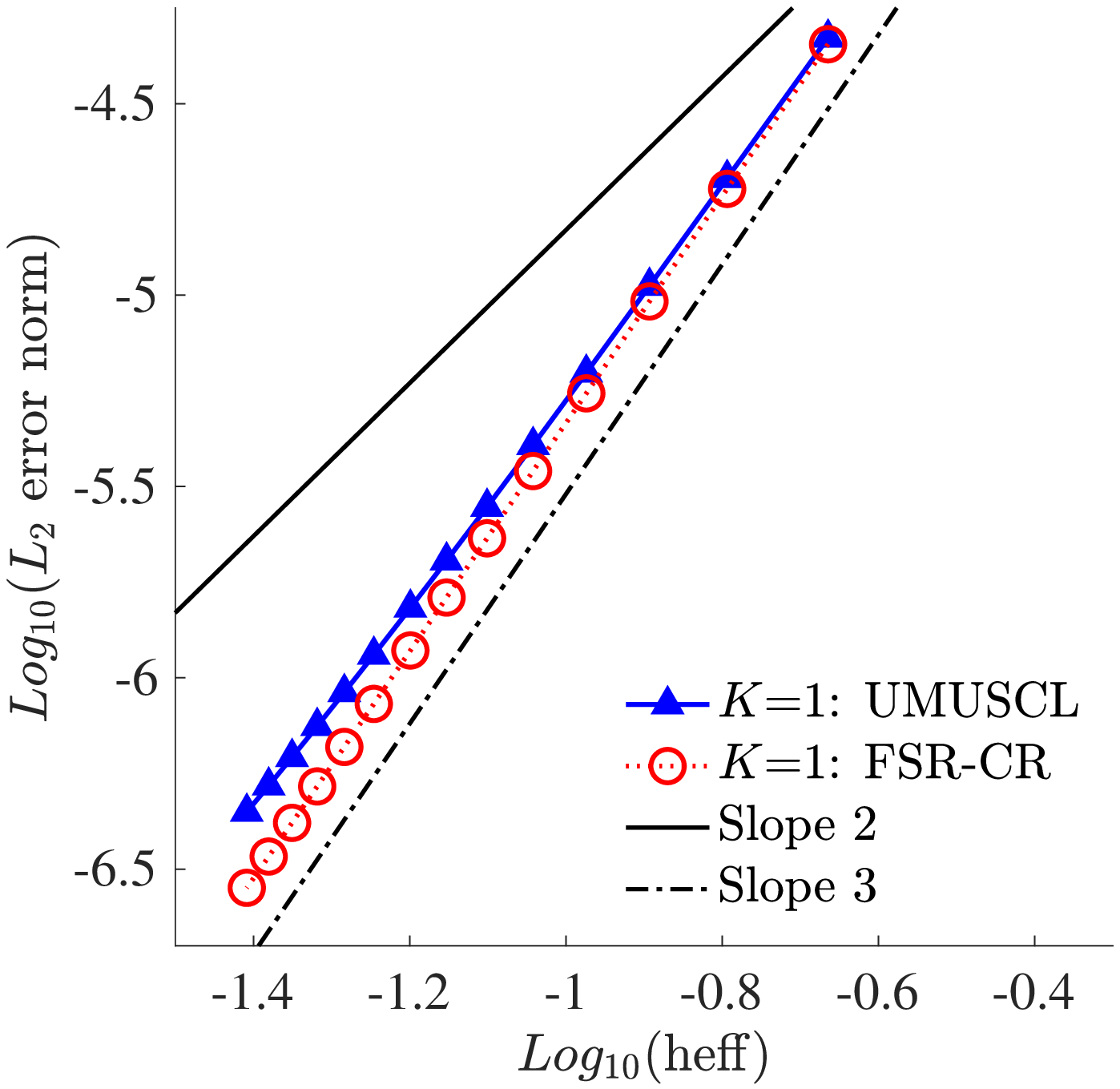}
          \caption{Error convergence: $K=1$.}
          \label{fig:inv_uns_err_mms_error_K1}
      \end{subfigure}
      \hfill
      \begin{subfigure}[t]{0.32\textwidth}
        \includegraphics[width=\textwidth,trim=0 0 0 0,clip]{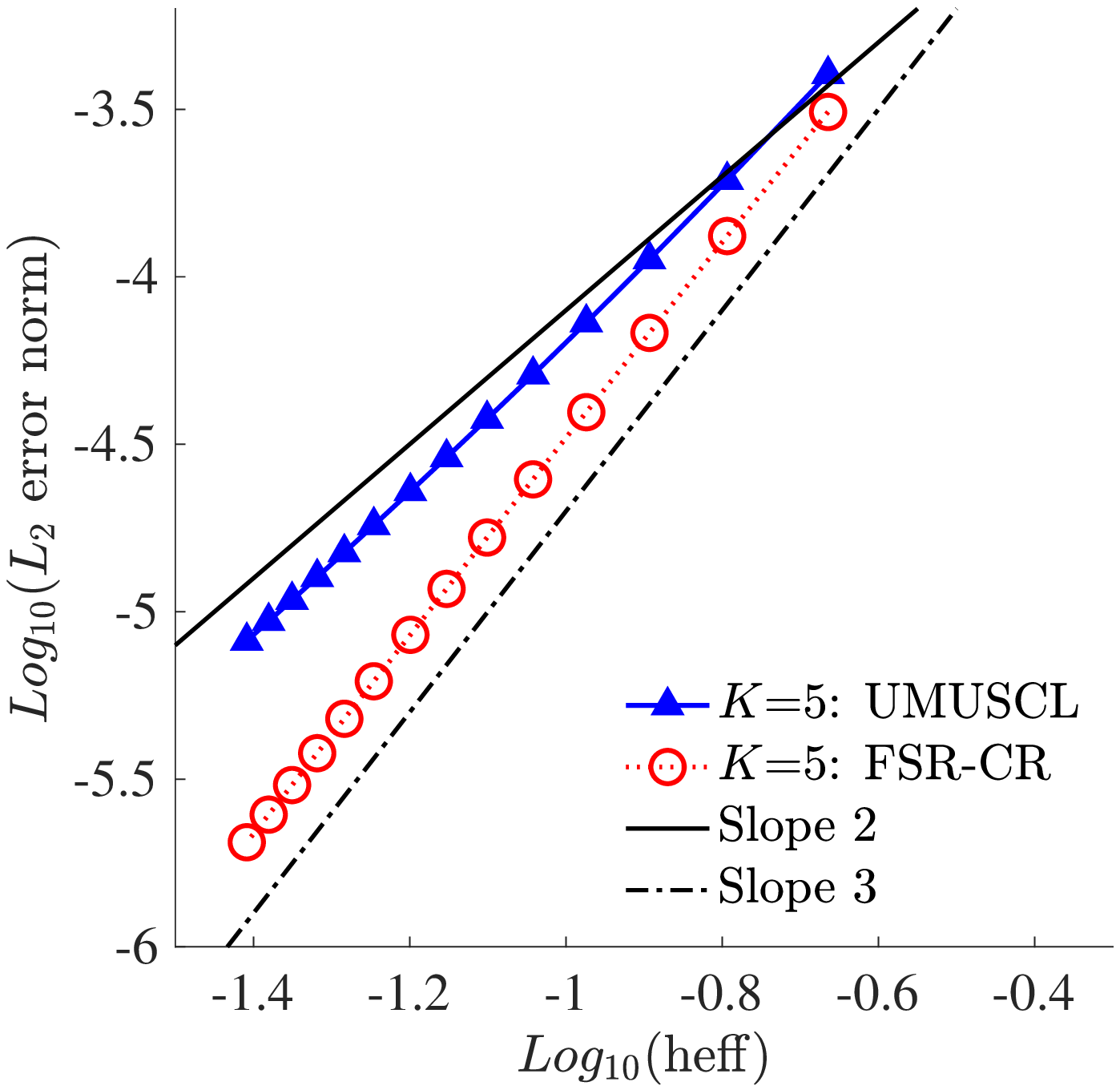}
          \caption{Error convergence: $K=5$.}
          \label{fig:inv_uns_err_mms_error_K5}
      \end{subfigure}
      \caption{Error convergence study for the unsteady Euler equations in two dimensions.}
\label{fig:inv_uns_err_mms} 
\end{figure}

Finally, we consider an inviscid-vortex problem frequently used for accuracy verification
\cite{Burg_etal:AIAA2003-3983,yang_harris:AIAAJ2016,fun3d_website,yang_harris:CCP2018,Barakos:IJNMF2018,ZhongSheng:CF2020}, which
is an exact solution to the Euler equations without forcing terms. Following 
Refs.\cite{yang_harris:AIAAJ2016,fun3d_website,yang_harris:CCP2018}, we set the exact solutions as 
\begin{eqnarray}
u =  u_\infty  -  \frac{ K \overline{y} }{2 \pi }   \exp \left(  \frac{1-\overline{r} ^2 }{2} \right) ,  \quad
v  = v_\infty   + \frac{ K \overline{x} }{2 \pi }   \exp \left(  \frac{1-\overline{r} ^2 }{2} \right),
\end{eqnarray}
and
\begin{eqnarray}
T =  1  -  \frac{ K^2 (\gamma-1)  }{8 \pi^2 }   \exp \left(  1-\overline{r} ^2   \right) , \quad
\rho = T^{  \frac{1}{\gamma-1}  } , \quad
p = \frac{ \rho ^{  \gamma   }  }{\gamma}    ,
\end{eqnarray}
where $\overline{x} = x  - u_\infty t$, $\overline{y} = y  - v_\infty t$, $\overline{r} ^2 =\overline{x}^2 + \overline{y}^2 $, and $(u_\infty,v_\infty)=(0.2,0.0)$.  
The initial solution at $t=0$ is shown in Figure \ref{fig:inv_uns_err_mms_gridsol}. Here, the parameter $K$ corresponds to $\epsilon$ in the previous cases; 
we would expect that the unexpected linearization 
occurs when $K$ is small. Two cases are considered: (1)$K=1$ as in Refs.\cite{yang_harris:AIAAJ2016,fun3d_website,yang_harris:CCP2018} (the one in Ref.\cite{ZhongSheng:CF2020} is much smaller) and (2)$K=5$ as in Ref.\cite{Barakos:IJNMF2018}. 
As we will see, $K=1$ is small enough for the unexpected linearization to occur, and $K=5$ is large enough to reveal second-order accuracy of the UMUSCL scheme.


For our purpose, it suffices to perform the calculation for a short time. Thus, 
we compute the solution at the final time $t_f=1.0$ with the three-stage  
SSP Runge-Kutta scheme \cite{SSP:SIAMReview2001} for the total of 1000 time steps with a constant time step $\Delta t = 0.001$,
which is so small that errors are dominated by the spatial discretization. To verify the spatial order of accuracy, we perform the computation
over a series of $n$$\times$$n$ regular quadrilateral grids, where $n=48, 64, 80, 96,112, 128, 144, 160, 176, 192, 208, 224, 240, 256$. 
The coarsest grid is shown with pressure contours in Figure \ref{fig:inv_uns_err_mms_gridsol}.
 
Error convergence results for $K=1$ are shown in Figure \ref{fig:inv_uns_err_mms_error_K1}. 
Here, for consistency with the results reported in Refs.\cite{yang_harris:AIAAJ2016,fun3d_website,yang_harris:CCP2018}, we compute the $L_2$ error norm for the pressure. 
As can be seen, the UMUSCL scheme exhibits nearly third-order accuracy with a slight deterioration on fine grids. Also plotted is error convergence
obtained with the FSR-CR scheme. As can be clearly seen, it is genuinely third-order with no sign of deterioration. 
Figure \ref{fig:inv_uns_err_mms_error_K5} shows the results obtained for $K=5$. 
It is clearly seen now that  the UMUSCL scheme quickly deteriorates to second-order. It is not genuinely third-order as expected. 
On the other hand, the FSR-CR scheme maintains third-order accuracy.

These results indicate that high-order error convergence reported in Refs.\cite{yang_harris:AIAAJ2016,yang_harris:CCP2018} 
is not genuine, but also that the UMUSCL scheme of Refs.\cite{yang_harris:AIAAJ2016,yang_harris:CCP2018} can be made genuinely 
high-order by flux reconstruction, which can be performed efficiently with the chain rule as suggested.

%

\section{Conclusions}
\label{conclusions}
  
We have shown that the UMUSCL scheme of Burg \cite{burg_umuscl:AIAA2005-4999}, which is 
defined with point-valued numerical solutions, is second-order accurate at best for nonlinear equations 
if used with point-valued time derivatives and source/forcing terms as in Refs.\cite{yang_harris:AIAAJ2016,yang_harris:CCP2018,DementRuffin:aiaa2018-1305}. 
As shown, the UMUSCL scheme is equivalent to the high-order conservative finite-difference scheme of Shu and Osher \cite{Shu_Osher_Efficient_ENO_II_JCP1989} for linear equations, but not for nonlinear equations. Thus, it cannot be high-order for nonlinear equations unless the flux is directly reconstructed. 
Nevertheless, as we have shown, it can exhibit third-order error convergence for a manufactured solution with a small perturbation
because a target nonlinear equation is effectively linearized by such a solution. An estimate of a critical value of the perturbation parameter has been
derived for the Burgers equation and its validity has been confirmed by numerical experiments. Similar results have been obtained for the Euler equations
in one dimension, which suggests the perturbation parameter for the velocity should be of the same order as those for other variables
to avoid the {\color{black} false}  high-order error convergence. Finally, we have demonstrated {\color{black} false}  third-order accuracy of UMUSCL
and genuine third-order accuracy of FSR (a flux and solution reconstruction scheme) for a steady problem and an unsteady inviscid vortex transport problem widely used for accuracy 
verification. In conclusion, any point-wise numerical scheme based on fluxes evaluated with reconstructed solutions cannot be high-order for nonlinear equations 
and accuracy verification must be carefully performed with exact solutions of a significant variation in order to avoid {\color{black} false}  high-order error convergence.

Thus far, we have clarified third-order accuracy of the MUSCL scheme \cite{Nishikawa_3rdMUSCL:2020}, third-order accuracy of the 
QUICK scheme \cite{Nishikawa_3rdQUICK:2020}, and the {\color{black} false}  third-order accuracy of the UMUSCL scheme; we are now ready to clarify economical high-order 
unstructured-grid schemes and identify one of the most efficient third- and higher-order schemes for practical simulations, which will be discussed in detail in a subsequent paper.

\addcontentsline{toc}{section}{Acknowledgments}
\section*{Acknowledgments}

The author gratefully acknowledges support from Software CRADLE, part of Hexagon, the U.S. Army Research Office 
under the contract/grant number W911NF-19-1-0429 with Dr. Matthew Munson as the program manager, and 
 the Hypersonic Technology Project, through the Hypersonic Airbreathing Propulsion Branch of the NASA Langley
 Research Center, under Contract No. 80LARC17C0004.

\addcontentsline{toc}{section}{References}
\bibliography{../../../bibtex_nishikawa_database}
\bibliographystyle{aiaa}

 
\end{document}